\documentclass[11pt]{amsart}
\usepackage{amssymb}
 \newcommand{\Z}{{\mathbb Z}}
 
 \newcommand{\F}{{\mathbb F}}
 \newcommand{\Q}{{\mathbb Q}}
 \newcommand{\R}{{\mathbb R}}
 \newcommand{\K}{{\mathbb K}}
 \newtheorem{theorem}{Theorem}[section]
 \newtheorem{lemma}[theorem]{Lemma}
 \newtheorem{proposition}[theorem]{Proposition}
 \newtheorem{corollary}[theorem]{Corollary}
 
 \newtheorem{example}[theorem]{Example}

\newcommand{\smod}{\hspace*{-0.2cm}\mod}
\newcommand{\pdk}[1]{$P\hspace{-0.05cm}D_{\K}(#1)$}
\newcommand{\pdq}[1]{$P\hspace{-0.05cm}D_{\Q}(#1)$}
\newcommand{\pdfp}[1]{$P\hspace{-0.05cm}D_{\F_p}(#1)$}
\newcommand{\spd}[2]{$P_{#1}(#2)$}
\newcommand{\sspd}[2]{$S\hspace{-0.05cm}P_{#1}(#2)$}
\newcommand{\swpd}[2]{$W\hspace{-0.1cm}P_{#1}(#2)$}

\renewcommand{\P}{Poincar\'e\ }
\newcommand{\oi}[1]{\Omega_{inv}^{{#1}}(M)}
 \def\Box  {\hfill \thinspace\vbox{\hrule height .5pt \hbox{\vrule
   width .5pt \vbox to 7pt{\hbox to 3.5pt{}} \vrule width .5pt}
   \hrule height 0pt depth .5pt}}

\title{Torus and $\Z/p$ actions on manifolds}

\author{Adam S. Sikora}

\date{}

\begin{document}

\maketitle
\begin{abstract} Let $G$ be either a finite cyclic group of
prime order or $S^1.$ We show that if $G$ acts on a manifold or,
more generally, on a \P duality space $M$, then each term of the
Leray spectral sequence of the map $M\times_G EG\to BG$ satisfies a properly
defined ``\P duality.'' As a consequence of this fact we obtain
new results relating the cohomology groups of $M$ and $M^G.$
We apply our results to study group actions on $3$-manifolds.
\end{abstract}

\subsection{Introduction}
\label{intro} If $G$ is either a finite cyclic group of prime
order, $\Z/p,$ or $S^1$ acting on a space $M$ then the
$G$-equivariant cohomology of $M$ can be calculated from the Leray
spectral sequence of the map $M\times_G EG\to BG.$ If $G=S^1$ and
$M$ is a  \P duality space then the components of the second term
of this spectral sequence (the cohomology groups of $M$) satisfy
\P duality. We show that if $M^{S^1}\ne \emptyset$ then each term
of this spectral sequence satisfies a properly defined ``\P
duality.''  Similarly, all terms of the corresponding Leray-Serre
and Cartan spectral sequences satisfy \P duality. These statements
and similar statements for $\Z/p$-actions will be formulated
precisely in Sections \ref{sPD-S^1} and \ref{sPD-Zp}.

Using the notion of \P duality for a spectral sequence we prove
new results relating the cohomology of any \P duality space with a
torus or $\Z/p$ action, to the cohomology of the fixed point set
of this action, see Theorems \ref{torus}, \ref{Zp}, \ref{Zp-Fp}.

Since this work was motivated by a conjecture concerning group
actions on $3$-manifolds, we devote Sections
\ref{sS^1-3}-\ref{sZp-3} to discuss the consequences and the
ramifications of the above results in $3$-dimensional topology.

{\bf Acknowledgments:}
We would like to thank C. Allday, A. Adem, A. Giacobbe, B. Hanke, V. Puppe,
J. A. Schafer, M. Sokolov and the referee of this paper
for many helpful comments.
For further comments and a different approach to our results
see \cite{AHP}, which was inspired by this paper.

For the reader who is unfamiliar with group cohomology,
equivariant cohomology, spectral sequences, or basic facts about
group actions we suggest
\cite{Adem,AP,Bredon-groups,Brown,McCleary, Weibel} as good
sources of information on these subjects.

Throughout this paper we will consider paracompact spaces $X$ of
finite cohomological dimension (over $\Z$) only. If $X$ is a
manifold or CW-complex then the cohomological dimension of $X,$
$cd\, X,$ is equal to $dim\, X.$ For more information on $cd\, X$
see \cite{Bredon-sheaf}. By $H^*$ we will denote the sheaf
cohomology groups with supports in closed sets. Recall that sheaf
cohomology theory with constant coefficients agrees with
Alexander-Spanier and \v Cech cohomology for paracompact spaces.
Let $b^i(X)=dim_\Q\, H^i(X;\Q).$

We say that a connected topological space $X$ is a \pdk{n}-space
(\P duality space of formal dimension $n$ with respect to
coefficients in a field $\K$) if $H^i(X;\K)=0$ for $i>n,$
$H^n(X;\K)=\K,$ and for all $0\leq i\leq n$ the cup product
$$H^i(X;\K)\times H^{n-i}(X;\K)\stackrel{\cup}{\longrightarrow}
H^n(X;\K)\cong \K$$ is a non-degenerate bilinear form. We also assume that
$dim_{\K}\, H^*(X;\K)<\infty.$ \pdk{n}-spaces will be usually denoted
by letter $M.$

\subsection{Torus actions}
\label{sTorus} We say that a torus $T$ action on a topological
space $X$ has finitely many connective orbit types (FMCOT) if the
set $\{(T_x)^0: x\in X\}$ is finite. Here $T_x$ denotes the
stabilizer of $x,$ $\{t\in T: tx=x\},$ and $(T_x)^0$ denotes the
connected component of identity of $T_x.$ Note that each
$S^1$-action has FMCOT.

\begin{theorem}\label{torus}
If a torus $T$ action on a \pdq{n}-space $M$ has FMCOT and
\begin{equation}
\begin{minipage}{4in}
either $n$ is even or\\
$M^{T}\ne \emptyset$ and $b^i(M)=0$ for all even $i,$
$0<i\leq \frac{1}{2}(n-1)$
\end{minipage}
\label{S^1-cond}
\end{equation}
then
$$\sum_i b^i(M^T)\equiv \sum_i b^i(M)\quad \mod 4.$$
\end{theorem}

It will be seen in Sections \ref{sS^1-3} and \ref{smore_ex}
that condition (\ref{S^1-cond}) is necessary.
The proof of Theorem \ref{torus} is given in Section \ref{sTorus-proof}.

The following well-known formulas provide additional information relating
$X$ to $X^T$ for any torus action with FMCOT:\\
\begin{equation}
\chi(X^T)= \chi(X),
\label{euler-torus}
\end{equation}
\begin{equation}
\sum_{i=0}^{\infty} b^{k+2i}(X^T)\leq \sum_{i=0}^{\infty}
b^{k+2i}(X), \label{inequality-torus}
\end{equation}
for all $k,$ cf. Theorems 3.1.13 and 3.1.14 in \cite{AP}. Here,
$\chi$ denotes the Euler characteristic.


\subsection{$\Z/p$-actions}
\label{sZp}


Let $\F_p$ denote the field of $p$ elements. We are going to see
that if $\Z/p$ acts on $X$ then the numbers
$$t^i(X)=dim_{\F_p} H^2(\Z/p,H^i(X;\F_p))$$
play similar role in relating $X$ to $X^{\Z/p}$ as the Betti numbers
in the study of torus actions.
For example, if $X$ is a finite dimensional
$\Z/p$-CW complex or a finitistic space then
(\ref{inequality-torus}) corresponds to
\begin{equation}
\sum_{i=0}^{\infty} t^{k+i}(X^{\Z/p})\leq \sum_{i=0}^{\infty}
t^{k+i}(X),
\label{inequality_of_t's}
\end{equation}
which holds for all $k,$ see \cite[Corollary
4.6.16]{AP}\footnote{In order to deduce (\ref{inequality_of_t's})
from \cite[Corollary 4.6.16]{AP} we need to notice that if $\Z/p$
acts on an $\F_p$-vector space $N$ then all Tate cohomology groups
$\hat{H}^i(\Z/p,N)$ are equal to $H^2(\Z/p,N).$ This can be proved
using Herbrand quotient or using the classification of
$\F_p[\Z/p]$-modules given in the next section.}.

Unlike for $S^1$-actions, the induced $\Z/p$-action on the
cohomology groups of $X$ may be non-trivial. For that reason, the
results for $\Z/p$-actions analogous to (\ref{euler-torus}) and
Theorem \ref{torus} can be formulated and proved only if the
action of $\Z/p$ on $H^*(X;\F_p)$ is {\em nice}.\vspace*{.1in}

\noindent{\bf Definition} {\em An action of $\Z/p$ on an
$\F_p$-vector space $N$ is nice if $N$ decomposes as
$\F_p[\Z/p]$-module into $T\oplus F,$ where $T$ and $F$ are
trivial and free $\F_p[\Z/p]$-modules respectively, ie.
$T=\bigoplus \F_p,$ $F=\bigoplus \F_p[\Z/p].$ (In particular,
trivial actions are nice.) We say that $\Z/p$ acts nicely on $X$
if the induced $\Z/p$-action on $H^*(X;\F_p)$ is nice. Note that
if $H^*(X;\F_p)=T^*\oplus F^*$ then
\begin{equation}\label{t_for_trivial}
t^i(X)=dim_{\F_p}\, T^i.
\end{equation}}

\begin{proposition}\label{p-euler-zp}(Proof in Section
\ref{speuler-zp}) If $p\ne 2$ and $\Z/p$ acts nicely on a space
$X$ such that $H^*(X;\Z)$ has no $p$-torsion then the following
version of the Euler characteristic formula holds:
$$ 
\chi_t(X^{\Z/p})= \chi_t(X), 
$$ 
where $\chi_t(X)=\sum_i (-1)^i\, t^i(X).$
\end{proposition}

For completeness, we recall also the classical formula, (see
\cite[Theorem III.4.3]{Bredon-groups}):
$$\chi(X)-\chi(X^{\Z/p})=p(\chi(X/\Z/p)-\chi(X^{\Z/p})),$$
which holds if $X$ is a finite dimensional
$\Z/p$-CW complex or a finitistic space.

If $\Z/p$ acts on a \pdfp{n}-space $M$ then
$t^i(M)=t^{n-i}(M)$ by Corollary \ref{dual_mod}.
Moreover, we have the following
counterpart of Theorem \ref{torus} for $\Z/p$-actions.

\begin{theorem}\label{Zp}(Proof in Section \ref{sZp-Z-proof})
Let $\Z/p$ act nicely on a $PD_{\F_p}(n)$-space $M$ with no
$p$-torsion in $H^*(M;\Z).$ If $p\ne 2,$ and
\begin{equation}
\begin{minipage}{4in}
either $n$ is even or\\
$M^{\Z/p}\ne \emptyset$ and $t^l(M)=0$ for all even
$0<l\leq \frac{1}{2}(n-1),$
\end{minipage}
\label{Zp-cond}
\end{equation}
then
$$\sum_i t^i(M^{\Z/p})\equiv \sum_i t^i(M)\quad \mod 4.$$
\end{theorem}

The proof of the above theorem for even $n$ is based on
Proposition \ref{p-euler-zp} and the standard properties of
\P duality spaces. For odd $n,$ the proof uses the notion
of \P duality of spectral sequences (defined in Section
\ref{sPD}) applied to the Leray spectral sequence associated with
the $\Z/p$-action on $M.$

The next result shows the assumption about the lack of $p$-torsion
in $H^*(M;\Z)$ can be replaced by the following condition:

\begin{equation}\label{cond}
\begin{minipage}{4in}
$\Z/p$ acts on a \pdfp{n}-space $M$ such that $d_r^{kl}=0$ for all
odd $r>1$ and all $k\geq n$ in the Leray spectral sequence of the
map $M\times_{\Z/p} E\Z/p\to B\Z/p$ with coefficients in $\F_p.$
\end{minipage}
\end{equation}

\begin{theorem}\label{Zp-Fp} (Proof in Section \ref{sZp-Fp-proof})
Let $p\ne 2$ and let $\Z/p$ act nicely on a \pdfp{n}-space $M$ in
such a way that conditions (\ref{Zp-cond}) and (\ref{cond}) hold.
If $M^{\Z/p}\ne \emptyset$ then
$$\sum_i t^i(M^{\Z/p})\equiv \sum_i t^i(M)\quad \mod 4.$$
\end{theorem}

We will see in Proposition \ref{SPD_3mflds} that all nice
$\Z/p$-actions on \pdfp{n}-spaces for $n\leq 3$ satisfy
(\ref{cond}). Furthermore, one can show that if $p\ne 2$ and
$\Z/p$ acts nicely on a \pdfp{n}-space $M$ with $M^{\Z/p}\ne
\emptyset$ then $d_r^{kl}=0$ for all odd $r>1$ and all odd $k\geq
n.$ Motivated by the above results, we conjectured in the previous
version of this paper that condition (\ref{cond}) holds for all
nice $\Z/p$-actions on \pdfp{n}-spaces with a non-empty fixed
point set for $p\ne 2.$ Recently B. Hanke showed that although
this conjecture is not true in general, it does hold under the
additional assumption that $H^*(M,\Z)$ does not contain $\Z/p$ as
a direct summand, \cite{Hanke}.

All other assumptions of Theorem \ref{Zp} are necessary. Examples
given in Sections \ref{sZp-3} and \ref{smore_ex} show that
condition (\ref{Zp-cond}) cannot be dropped. We will also see that
Theorem \ref{Zp} fails if the $\Z/p$-action on $M$ is not nice and
$n>2.$ However, a much stronger statement holds for
$2$-dimensional manifolds\vspace*{.1in}

\noindent{\bf Theorem} \cite{Bryan} {\em
If $\Z/p$ acts on a connected surface $F,$ $F^{\Z/p}\ne \emptyset,$
then this action has $2+dim_{\F_p}H^1(\Z/p, H_1(F,\F_p))$ fixed
points.}


\subsection{$S^1$-actions on $3$-manifolds}
\label{sS^1-3}


Since this paper was motivated by a conjecture concerning group
actions on $3$-manifolds, we devote this and the next subsection
to present consequences and ramifications of our results to such
actions.

By the slice theorem, if $G=S^1$ or $\Z/p$ acts smoothly on a
closed, oriented, smooth manifold $M$ with fixed points then
$M^G$ is a disjoint union of closed, orientable
submanifolds of even codimension. Therefore, if $dim\, M=3$ then
$M^G$ is a union of embedded circles.

By Theorem \ref{torus} and by (\ref{inequality-torus}) we have

\begin{corollary}\label{S^1-3}
If $S^1$ acts smoothly on a connected, closed, orientable $3$-manifold $M,$
$M^{S^1}\ne \emptyset,$ then
$M^{S^1}$ is a union of $s$ circles, where
\begin{itemize}
\item $s\leq 1+b_1(M)$ and \item $s\equiv 1+b_1(M)$ mod $2.$
\end{itemize}
\end{corollary}

This corollary can be also deduced from the classification of
$S^1$-actions on $3$-manifolds, \cite{Ra, OR}.

The statement of the corollary cannot be improved. Namely, given
two integers $s,b$ such that $0< s\leq 1+b$ and $s\equiv 1+b$ mod
$2,$ there is a $3$-manifold $M$ with an $S^1$-action such that
$b_1(M)=b$ and $M^{S^1}$ is a union of $s$ circles. Such a
manifold can be constructed as follows. Let $F_{g,s}$ denote a
surface of genus $g=\frac{b+1-s}{2}$ with $s$ boundary components,
and let $M_0=F_{g,s}\times S^1.$ The boundary of $M_0$ is a union
of $s$ tori and $b_1(M_0)=b_1(F_{g,s})+b_1(S^1)=b+1.$ Choose $s$
points $p_1,...,p_s\in
\partial F_{g,s},$ each lying in a different component of $\partial
F_{g,s}.$ Now, attach $s$ solid tori to $M_0$ along their boundaries, in a
such a way that the meridian of the $i$-th
solid torus is identified with $p_i\times S^1.$ We denote the closed
manifold obtained in this way by $M.$ Note that after attaching the first
solid torus the first
Betti number decreases by $1,$ but after attaching the next tori, it stays
unchanged. Therefore $b_1(M)=b.$ Obviously, the $S^1$-action on
$M_0$ extends on $M$ and the fixed point set of the action
is composed of the cores of the solid tori. Hence, $M^{S^1}$ has exactly
$s$ components.


\subsection{$\Z/p$-actions on $3$-manifolds}
\label{sZp-3}


If $\Z/p$ acts on a closed, connected, orientable $3$-manifold $M$ then
(\ref{inequality_of_t's}) for $k=0$ implies that $M^{\Z/p}$ is a union of
at most $1+t^1(M)$ circles.

Since by Proposition \ref{SPD_3mflds} $\Z/p$-actions on
$3$-manifolds with fixed points always satisfy condition
(\ref{cond}), the following result is a special case of Theorem
\ref{Zp-Fp}.

\begin{proposition}\label{Zp-3mfld}
Let $p\ne 2.$ If $\Z/p$ acts nicely on a closed, connected,
orientable $3$-manifold $M$ and if $M^{\Z/p}$ is composed of $s$
circles, $s\ne 0,$ then $s\equiv 1+t^1(M)$ mod $2.$
\end{proposition}

The above proposition answers in affirmative a conjecture of M.
Sokolov concerning $p$-periodic $3$-manifolds, i.e. manifolds with
a $\Z/p$-action with exactly one circle of fixed points. M.
Sokolov conjectured the following statement.

\begin{proposition}
If $M$ is $p$-periodic then $H^1(M;\F_p)\ne \F_p.$
\end{proposition}

\begin{proof}
If $H_1(M;\F_p)$ is a $1$-dimensional vector space then the
$\Z/p$-action on $H_1(M;\F_p)$ is trivial, and therefore $t^1(M)=1.$
Hence, $M^{\Z/p}=\emptyset$ or $S^1\cup S^1.$
\end{proof}

The statement of the above proposition makes an impression that it
could be easily proved by elementary means of algebraic topology,
Smith theory, or $3$-dimensional topology. We do not know any
short proof of it, and we encourage the reader to try to find one
by himself, in order to realize that this is not easy. After we
proved the above proposition, J. H. Przytycki and M. Sokolov
\cite{PS} found a different proof of it, which avoids using
equivariant cohomology at the expense of an elaborate application
of surgery theory.

Now we are about to show that Proposition \ref{Zp-3mfld} does not
hold if we drop any of its assumptions. Observe first that any
free $\Z/p$-action on $S^3$ is nice and $s=0$ and $t^1(S^3)=0$ for
such action. Therefore the assumption $s\ne 0$ in Proposition
\ref{Zp-3mfld} is necessary.

If $p=2$ or if the action is not nice then the conclusion of
Proposition \ref{Zp-3mfld} fails as well. To see this, we need to
understand the possible $\Z/p$-actions on $H_1(M,\F_p).$ If $\Z/p$
acts on $M$ then $H_1(M,\F_p)$ considered as a module over
$R=\F_p[\Z/p]$ decomposes as a direct sum of indecomposable
$R$-modules. We will see in Section \ref{sreps} that each
indecomposable $R$-module is isomorphic to
$$V_i=R/(t-1)^i=\F_p[t]/(t-1)^i,$$ for a unique $i$ between $1$
and $p.$ (Here, $R=\F_p[t]/(t^p-1)= \F_p[t]/(t-1)^p.$) A
$\Z/p$-action on an $\F_p$-vector space $N$ is nice if $N$
decomposes as an $R$-module into a sum of $V_1$'s and $V_p$'s. By
Corollary \ref{dual_mod}, a $\Z/p$-action on a $3$-manifold $M$ is
nice if and only if the induced $\Z/p$-action on $H_1(M,\F_p)$ is
nice.

A $3$-manifold with a $\Z/p$-action which is not nice can be
constructed as follows. Let $S^3_1$ and $S^3_2$ be two $3$-spheres
with some (not necessarily the same) $\Z/p$-actions. Choose
$3$-balls $B_1\subset S_1^3,$ $B_2\subset S_2^3,$ such that the
orbit of $B_i,$ $\bigcup_{g\in \Z/p} gB_i,$ for $i=1,2,$ is
composed of $p$ disjoint balls (on which $\Z/p$ acts freely). For
all $g\in \Z/p,$ remove the interiors of the balls $gB_1,$ $gB_2,$
from $S^3_1$ and $S^3_2$ respectively. Next, choose an arbitrary
homeomorphism $\Psi:\partial B_1 \to \partial B_2$ and identify
$g\partial B_1$ with $g\partial B_2,$ for any $g\in \Z/p,$ via
$g\Psi g^{-1}.$ This construction gives a closed, orientable
$3$-manifold $M_p,$ with the cyclic group, $\Z/p,$ acting on it.

The proofs of the following remarks are easy and left to the
reader.\\

\noindent{\bf Remarks}
\begin{enumerate}
\item[(i)] $M_p\simeq \underbrace{(S^2\times S^1) \# ... \#
(S^2\times S^1)}_{p-1}.$
\item[(ii)] $H_1(M_p,\F_p)\cong V_{p-1}$
for any $\Z/p$-action on $M_p$ constructed as above.
In particular, none of the $\Z/p$-actions on $M$ is nice for $p\ne
2.$
 \item[(iii)] On the other hand, since the only indecomposable
$\F_2[\Z_2]$-modules are the trivial module, $V_1,$ and the free
module, $V_2,$ all $\Z_2$-actions on vector spaces over $\F_2$ are
nice. In particular, any $\Z_2$-action on $M_2$ is nice.
\item[(iv)] Since $\Z/p$ can act on $S^3$ with $(S^3)^{\Z/p}= S^1$
or $\emptyset,$ there exist $\Z/p$-actions on $M_p$ with
$M_p^{\Z/p}=\emptyset, S^1,$ and $S^1\cup S^1.$ \item[(v)] By (ii)
and (iv) the statement of Proposition \ref{Zp-3mfld} fails for
$\Z/p$-actions which are not nice. \item[(vi)] By (ii) and (iii)
the statement of Proposition \ref{Zp-3mfld} fails for $p=2.$
\end{enumerate}

\subsection{More examples}
\label{smore_ex}

Theorems \ref{torus} and \ref{Zp} may be useful for
studying group actions on products of spheres. On the other hand,
the analysis of examples of such actions shows that all the
assumptions of Theorems \ref{torus} and \ref{Zp} are necessary.

\begin{example}\label{ex1}
Bredon in \cite[VII \S 10]{Bredon-groups} constructs a circle action on
$M=S^3\times S^5 \times S^9$ with the fixed point set $M^{S^1}$ being an
$S^7$-bundle over $S^3\times S^5$ with $b_i(M^{S^1})=1$ for
$i=0,3,5,10,12,15$ and $b_i(M^{S^1})=0$ for all other $i.$
This circle action does not satisfy the conclusion of Theorem \ref{torus}.
Therefore condition (\ref{S^1-cond}) is necessary.
\end{example}

\begin{example} There is a $\Z_3$-action on $M=S^n\times S^n,$ for
$n=1, 3$ or $7,$ which is not nice and for which $M^{\Z_3}=
(point + S^{n-1}),$ see \cite[VII \S 9]{Bredon-groups}. This shows that
the restriction in Theorem \ref{Zp} to nice actions is necessary.
The action can be constructed as follows: Let $R$ be the ring of
complex numbers, quaternions, or Cayley numbers for $n=1,3,7$ respectively.
Let $S$ denote the set of elements of norm $1$ in $R,$ $S\simeq S^n,$
and let $M$ be the space of all triples $(x,y,z)\in S\times S\times S$
such that $(xy)z=1.$ $M$ is homeomorphic to $S^n\times S^n$ and
since $$(xy)z=1 \Longleftrightarrow (yz)x=1 \Longleftrightarrow (zx)y=1$$
there is an action of $\Z_3$ on $M$ by cyclic permutations.
The fixed point set of this action
is $\{x\in R| x^3=1\}\simeq point + S^{n-1}.$
\end{example}

The notation $X\sim_p Y$ in the next example means that
$X$ and $Y$ are topological spaces with isomorphic cohomology rings with
coefficients in $\F_p.$

\begin{example} If $n\ne m$ and $n, m$ are both even or both odd,
or if the smaller of them is odd then any action of $\Z/p$ on $M=S^n\times
S^m$ is nice and condition
(\ref{Zp-cond}) is satisfied. In this situation Theorem
\ref{Zp} holds, and one can prove that $M^{\Z/p}$
is $\sim_p$-equivalent to one of the following spaces:
$S^q\times S^r, S^q+S^r, P^3(2q), (point+P^2(2q));$ see \cite[Thm. VII 9.1]
{Bredon-groups}.
Here $P^n(2q)$ denotes a space whose cohomology ring with coefficients
in $\F_p$ is $\F_p[x]/(x^{n+1}),$ and $deg\, x =2q.$
\end{example}

However, there are known examples of $X\sim_p S^n\times S^m$
which do not satisfy the assumptions of the example above
(i.e. $\min(n,m)$ is even and $\max(n,m)$ is odd)
and which admit a $\Z/p$-action with $X^{\Z/p}\sim_p S^q.$
Therefore (\ref{Zp-cond}) is a necessary condition
for Theorem \ref{Zp}.


\section{Classification of representations of $\Z/p$}
\label{sreps}


In this section we present a classification of all representations
of $\Z/p$ over $\F_p$ and over the ring of integers localized at
the prime ideal $(p),$ $\Z_{(p)}.$ This classification should help
the reader to better understand the possible $\Z/p$-actions on the
cohomology groups of $X.$ We will classify indecomposable modules
only, since all other modules are direct sums of these.

Note that $R=\F_p[\Z/p]$ is isomorphic to $\F_p[t]/(t^p-1),$ and
since $t^p-1=(t-1)^p$ mod $p,$ $R=\F_p[t]/(t-1)^p.$ Therefore,
$$V_i=R/(t-1)^i=\F_p[t]/(t-1)^i,$$ is an $R$-module for each
$i=1,...,p.$

\begin{proposition}\label{Zp-reps}
\begin{enumerate}
\item $V_1,...,V_p$ is the complete list of finitely generated
indecomposable $R$-modules. \item Each finitely generated
$R$-module, $N,$ decomposes as a finite sum
$$N=V_{i_1}\oplus ... \oplus V_{i_k},$$ where $1\leq i_1, ...,
i_k\leq p$ are unique up to a permutation.
\end{enumerate}
\end{proposition}

\begin{proof} (i) $\F_p[\Z/p]=\F_p[t]/(t-1)^p$
is a quotient of the ring of polynomials $\F_p[t],$ which is a
principal ideal domain. Every indecomposable module over
$\F_p[\Z/p]$ is also indecomposable over $\F_p[t]$ and, hence,
cyclic. Such modules
can be easily classified.\\
(ii) Note that $V_k/(t-1)^d=V_{min(k,d)}.$ Therefore,  the number
of components $V_d$ in $N$ is determined by difference in the
dimensions of the vectors spaces $N/(t-1)^d$ and $N/(t-1)^{d+1}.$
(Another way of proving the uniqueness of the decomposition of $N$
is by applying the Krull-Schmidt Theorem, \cite[14.5]{Curtis}).
\end{proof}

\begin{corollary}\label{dual_mod}
If $N$ and $N'$ are $\F_p[\Z/p]$-modules and $\Psi: N\times N'\to
\F_p$ is a non-degenerate bilinear $\Z/p$-equivariant form then
$N$ and $N'$ are isomorphic as modules.
\end{corollary}

\begin{proof}
$N'$ is the dual module to $N,$ ie. $N'$ is isomorphic to
$Hom_{\F_p}(N,\F_p),$ where $g\in \Z/p$ sends $f:N\to \F_p$ to the
homomorphism $x\to f(g^{-1}x).$ Since $(N_1\oplus N_2)'=N_1'\oplus
N_2',$ it is enough to assume that $N=V_k.$ The module $V_k'$ is
generated by the homomorphism given by $f(t^i)=\begin{cases} 1 &
\text{for $i=0$}\\ 0 & \text{for } 1\leq i\leq k-1.
\end{cases}$
Therefore $V_k'$ is a cyclic $R$-module and hence $V_k'=V_i$ for
certain $i.$ Now $i=k$ since $dim_{\F_p}\, V_k'=dim_{\F_p}\, V_k.$
\end{proof}

Note that a $\F_p[\Z/p]$-module $N$ is nice if it decomposes into
a sum of $V_1$'s and $V_p$'s.

Later we will need the classification of $\Z/p$-modules over
$\Z_{(p)}.$ By modifying the proof of Theorem 74.3 in
\cite{Curtis} one can show the following:

\begin{proposition}
Every indecomposable $\Z/p$-module over $\Z_{(p)}$ which is free
over $\Z_{(p)}$ is either
\begin{enumerate}
\item the trivial module, $\Z_{(p)},$ or \item the free module,
$\Z_{(p)}[\Z/p],$ or \item the ring of cyclotomic integers,
$\Z_{(p)}[\zeta_p]=\Z_{(p)}[x]/(1+x+...+x^{p-1}).$
\end{enumerate}

The action of the generator of $\Z/p$ on
$\Z_{(p)}[\zeta_p]$ is given by the multiplication by $\zeta_p.$
\end{proposition}

Note that $\Z_{(p)}\otimes \F_p=V_1,$ $\Z_{(p)}[\Z/p] \otimes
\F_p=V_p,$ and $\Z_{(p)}[\zeta_p]\otimes \F_p= V_{p-1}.$ Hence, if
$H^*(M;\Z)$ has no $p$-torsion then $H^*(M;\F_p)$ decomposes as a
sum of $V_1$'s, $V_{p-1}$'s, and $V_p$'s.

The following lemma will be needed in Section \ref{sPD-Zp}.

\begin{lemma}\label{nice}
If $p\ne 2$ and $\Z/p$ acts nicely on a space $X,$ with no
$p$-torsion in $H^*(X;\Z)$ then
$$H^k(\Z/p, H^*(X;\Z))=\begin{cases}H^2(\Z/p, H^*(X;\F_p)) &
\text{if $k$ is even}\\
0 & \text{if $k$ is odd,}\\
\end{cases}$$
for $k>0.$
\end{lemma}

\begin{proof}
By the universal coefficient theorem, $H^*(X;\Z_{(p)})$ is a free
$\Z_{(p)}$-module and $H^*(X;\F_p)=H^*(X;\Z_{(p)})\otimes \F_p.$
Therefore the sequence
$$0\to H^*(X;\Z_{(p)})\stackrel{\cdot p}{\to} H^*(X;\Z_{(p)})\to
H^*(X;\F_p)\to 0$$ is exact. By applying $H^*(\Z/p,\cdot)$ to that
sequence we get
$$0\to H^2(\Z/p,H^*(X;\Z_{(p)}))\to H^2(\Z/p,H^*(X;\F_p))\to
H^3(\Z/p,H^*(X;\Z_{(p)}))\to 0.$$ Now we use the classification of
$\Z/p$ modules over $\F_p$ and $\Z_{(p)}$ described above. Since
$\Z_{(p)}[\zeta_p]\otimes \F_p=V_{p-1},$ $H^*(X;\Z_{(p)})$ must be
a direct sum of trivial and free $\Z/p$-modules. Therefore,
$H^3(\Z/p, H^*(X;\Z_{(p)}))=0$ and, hence, $H^2(\Z/p,
H^*(X;\Z_{(p)}))= H^2(\Z/p, H^*(X;\F_p)).$ Since localization at
$(p)$ is an exact functor in the category of $\Z/p$-modules,
$H^2(\Z/p, H^*(X;\Z))=H^2(\Z/p, H^*(X;\Z_{(p)})).$ Finally, the
statement follows from the fact that $H^k(\Z/p, \cdot)$ is
$2$-periodic for $k>0.$
\end{proof}


\section{\P duality on spectral sequences}
\label{sPD}


Throughout this section we will make the following assumptions:
Let $\K$ be a field and $n$ be a positive integer. Let
$(E^{**}_*,d_*)$ be a spectral sequence whose each summand, $E^{pq}_r,$
for $r\geq 2,$ is a finite dimensional vector space
over $\K$ and $E^{kl}_2=0$ for $l<0$ and $l>n.$ Assume that
$E^{**}_r$ has a multiplicative structure for $r\geq 2,$ ie.
there is a graded commutative product on each term, $E^{**}_r,$ such
that $d_r^{**}$ is a derivation with respect to that product, and
the product on $E^{**}_{r+1}$ is induced from the product on $E^{**}_r.$
Additionally, assume the following condition about the $0$th row in
$E^{**}_*:$

\begin{minipage}{4.7in}
(ZR) $E^{*0}_2=E^{*0}_{\infty}.$ Equivalently, the differentials
$d^{*,r-1}_r: E^{*,r-1}_r\to E^{*+r,0}_r$ are $0$ for all $r\geq
2.$
\end{minipage}

\begin{proposition}\label{PD-main-prop}
Let $r\geq 2,$ $k,k'\in \Z$ and $l,l'\geq 0$ be such that
$E^{k+k',n}_r=E^{k+k',0}_r=\K,$ $l+l'=n,$ and the following
$\K$-bilinear maps are non-degenerate
\begin{equation}
E^{kl}_r\times E^{k',l'}_r\stackrel{\cdot}{\to} E^{k+k',n}_r=\K
\label{nondeg1}
\end{equation}

\begin{equation}
E^{k+r,l-r+1}_r\times E^{k'-r,l'+r-1}_r\stackrel{\cdot}{\to}
E^{k+k',n}_r=\K
\label{nondeg2}
\end{equation}

\begin{equation}
E^{k-r,l+r-1}_r\times E^{k'+r,l'-r+1}_r\stackrel{\cdot}{\to}
E^{k+k',n}_r=\K
\label{nondeg3}
\end{equation}
(A pairing $V\times W\to \K$ is non-degenerate if it induces
an isomorphism $V\to W^*.$ In particular, the pairing
$\{0\}\times \{0\}\stackrel{0}{\to} \K$ is non-degenerate.)

Assume additionally that (\ref{nondeg1}) and (\ref{nondeg2}) are
non-degenerate for $l=n,$ $l'=0.$ Under the above assumptions,
$E^{k+k',n}_{r+1}=\K$ and the pairing
\begin{equation}
E^{kl}_{r+1}\times E^{k',l'}_{r+1}\stackrel{\cdot}{\to} E^{k+k',n}_{r+1}=\K
\label{nondeg4}
\end{equation}
is non-degenerate.
\end{proposition}

\begin{proof}
Consider the diagram\\
\begin{equation}
\begin{tabular}{ccccc}
$E^{k-r,l+r-1}_r$ & $\stackrel{d^{k-r,l+r-1}_r}{\longrightarrow}$ &
$E^{kl}_r$ & $\stackrel{d^{kl}_r}{\longrightarrow}$ & $E^{k+r,l-r+1}_r$\\
$\Big \downarrow {\wr}$ & & $\Big \downarrow \wr$ & & $\Big \downarrow \wr$\\
$(E^{k'+r,l'-r+1}_r)^*$ \hspace*{-.1in} & $\stackrel{(d^{k'l'}_r)^*}
{\longrightarrow}$ \hspace*{-.1in} & $(E^{k'l'}_r)^*$ \hspace*{-.1in} &
$\stackrel{(d^{k'-r,l'+r-1}_r)^*}{\longrightarrow}$
\hspace*{-.1in}
& $(E^{k'-r,l'+r-1}_r)^*$\\
\end{tabular}
\label{diagram}
\end{equation}
in which the vertical isomorphisms are induced by the
non-degenerate pairings (\ref{nondeg1}), (\ref{nondeg2}), (\ref{nondeg3}),
$$E^{kl}_r\ni \alpha\quad \stackrel{\sim}{\longrightarrow} \quad
(x \to \alpha\cdot x)\in (E^{k'l'}_r)^*,$$
$$E^{k\pm r,l\mp r+1}_r\ni \beta \stackrel{\sim}{\longrightarrow}
(x \to \beta \cdot x)\in (E^{k'\mp r,l'\pm r-1}_r)^*.$$
$d_r^{kl}$ followed by the arrow pointing down
on the right side of the diagram sends $\alpha\in E^{kl}_r$ to the functional
$x\to d_r^{kl}(\alpha)\cdot x$ in $(E^{k'-r,l'+r-1}_r)^*,$ and the
arrow pointing down in the middle of the diagram followed by
$(d^{k'-r,l'+r-1}_r)^*$ sends $\alpha$ to
$x\to \alpha\cdot d_r^{k'-r,l'+r-1}(x).$
Since $n+r-1>n$ and $\alpha\cdot x\in E_r^{k+k'-r,n+r-1},$ $\alpha\cdot x=0$
and, hence, $d_r^{kl}(\alpha)\cdot x=\pm \alpha\cdot d_r^{k'-r,l'+r-1}(x).$
Therefore the right square in (\ref{diagram}) is commutative up to sign.
Similarly, we prove that the left square also commutes up to sign.
Since the top and the bottom rows are isomorphic chain
complexes, their cohomology is isomorphic, Hence, we have an
isomorphism
\begin{equation}
E^{kl}_{r+1}\ni \alpha\quad \longrightarrow
\quad (x \to \alpha\cdot x)\in (E^{k'l'}_{r+1})^*.
\label{prop-eq}
\end{equation}
By (ZR), $E^{k+k',0}_r=\K$ implies that $E^{k+k',0}_{r+1}=\K.$
Now (\ref{prop-eq}) for $l=n$ implies that $E^{k+k',n}_{r+1}=\K.$
Finally, by (\ref{prop-eq}), the pairing (\ref{nondeg4}) is non-degenerate.
\end{proof}

Now, we are going to use Proposition \ref{PD-main-prop} to prove the
existence of ``\P duality'' on spectral sequences.
This duality will be useful for the study of group actions on \P
duality spaces.

Let $(E^{**}_r,d_*)$ be a term of a spectral sequence with a
multiplicative structure such that $E^{*l}_r=0$ for all $l<0$ and
for all $l>n$ for a certain $n.$ We will consider $3$ types of
\P duality on $(E^{**}_r,d_r):$

 We say that $(E^{**}_r,d_r)$ satisfies {\bf \P
duality,} denoted by \spd{\K}{n}, if there exists $N$ such that
\begin{enumerate}
\item $E^{k*}_r=0$ for all odd $k>N$;
 \item $E^{k0}_r=E^{kn}_r=\K$ for all even $k> N$
 \item $E^{kl}_r\times E^{k',l'}_r\stackrel{\cdot}{\to} E^{k+k',n}_r=\K$ is
non-degenerate for all $l,l'\geq 0$ such that $l+l'=n$ and for all
even $k,k'\geq N.$
\end{enumerate}

We say that $(E^{**}_r,d_r)$ satisfies {\bf weak \P duality,}
denoted by \swpd{\K}{n}, if there exists $N$ such that

\begin{enumerate}
\item $E^{kn}_r=\K$ for all odd $k>N,$ and
 \item for all $k,k'>N$ of different
 parity and for all $0\leq l\leq n$ and $r\geq 2,$ the pairing $E^{kl}_r\times
E^{k',n-l}_r\stackrel{\cdot}{\to} E^{k+k',n}_r=\F_p$ is
non-degenerate.
\end{enumerate}

Finally, $(E^{**}_r,d_r)$ satisfies {\bf strong \P duality,}
\sspd{\K}{n}, if there exists $N$ such that
\begin{enumerate}
\item $E^{kn}_r=\K$ for all $k>N,$
 \item for all $k,k'>N$ such that at least one of them is even and for
all $0\leq l\leq n$ the pairing $E^{kl}_r\times
E^{k',n-l}_r\stackrel{\cdot}{\to} E^{k+k',n}_r=\K$ is
non-degenerate.
\end{enumerate}

The following statement follows by induction from Proposition
\ref{PD-main-prop}:

\begin{proposition}\label{PD-prop}
Let $(E^{**}_*,d_*)$ have a multiplicative structure and satisfy
(ZR).
\begin{enumerate}
 \item If $(E^{**}_r,d_r)$ satisfies \spd{\K}{n} for $r=2$ then
it satisfies \spd{\K}{n} for all $r> 2.$
 \item If $(E^{**}_r,d_r)$ satisfies \swpd{\K}{n} for $r=2$ then
it satisfies \swpd{\K}{n} for all $r> 2.$
 \item If $(E^{**}_r,d_r)$ satisfies \sspd{\K}{n} and $r$ is even
 then $(E^{**}_{r+1},d_{r+1})$ satisfies \sspd{\K}{n} as well.
\end{enumerate}
\end{proposition}


\begin{lemma}\label{prop_of_PD}
(i) If $r\geq 2$ and $(E^{**}_r,d_r)$ satisfies \spd{\K}{n} then
there exists $N,$ such that $E^{kl}_r\cong E^{k'l}_r,$ $rank\,
d_r^{kl}=rank\, d_r^{k'l},$ and $rank\, d_r^{kl}=rank\,
d_r^{k',n-l+r-1}$ for all even $k,k'>N$ and all $l.$\\
(ii) If
$r\geq 2$ and $(E^{**}_r,d_r)$ satisfies \sspd{\K}{n} then there
exists $N,$ such that $E^{kl}_r\cong E^{k'l}_r,$ for all $k,k'\geq
N.$ Additionally, if $r$ is even then $rank\, d^{kl}_r=rank\,
d^{k'l}_r,$ and $rank\, d^{kl}_r=rank\, d^{k',n-l+r-1}_r$ for all
$k,k'\geq n$ and $0\leq l\leq n.$
\end{lemma}

\begin{proof}(i)
$E^{kl}_r\cong (E^{k,n-l}_r)^*\cong E^{k'l}_r$ implies the first claim.
Since $d_r=0$ for $r$ odd, assume that $r$ is even. Note that the vertical
maps in (\ref{diagram}) are isomorphisms for $k,k',l,l'$ such that
$k,k'$ are even and sufficiently big and $l+l'=n.$
Hence
$rank\, d_r^{kl}=rank\, d_r^{k'-r,n-l+r-1}$ for any $0\leq l\leq n,$
and by substituting $k'$ for $k'-r$ we get
$$rank\, d_r^{kl}=rank\, d_r^{k',n-l+r-1}.$$
By applying this identity twice,
we get
$$rank\, d_r^{kl}=rank\, d_r^{k',n-l+r-1}=rank\, d_r^{k'l}.$$
The proof of (ii) is analogous.
\end{proof}


\subsection{\P duality for spectral sequences for $S^1$-actions}
\label{sPD-S^1}


\begin{proposition}{\bf PD for Leray spectral sequence}\label{PD-S^1}
If $S^1$ acts on a \pdq{n}-space $M$ and $M^{S^1}\ne \emptyset$
then the Leray spectral sequence of the map $\pi: M\times_{S^1}
ES^1\to BS^1$ with coefficients in $\Q$ satisfies condition
\spd{\Q}{n} for all $r\geq 2.$
\end{proposition}

\begin{proof} Since all cohomology groups are considered with
supports in closed sets, by \cite[Theorem IV.6.1]{Bredon-sheaf} we
have $E^{kl}_2= H^k(BS^1,{\mathcal H}^l(\pi;\Q)),$ where
${\mathcal H}^l(\pi;\Q)$ is the Leray sheaf of $\pi.$ Since $BS^1$
is simply connected, by the remark following Theorem IV.8.2 in
\cite{Bredon-sheaf}, ${\mathcal H}^l(\pi;\Q)$ is the constant
sheaf with stalks $H^l(X;\Q).$ Hence, $E^{kl}_2\cong
H^k(BS^1,\Q)\otimes H^l(M;\Q).$ The Leray spectral sequence of
$\pi,$ $(E^{**}_*,d_*),$ has a multiplicative structure -- see eg.
\cite[IV.6.5]{Bredon-sheaf} or \cite[XII \S 3.2]{Swan}. An
argument similar to that used in the proof of Theorem III.15.11 in
\cite{Bott-Tu} shows that $E^{**}_2$ and $H^*(BS^1,\Q)\otimes
H^*(M;\Q)$ are isomorphic as algebras. Since $H^*(BS^1)=\Q[t],$
where $deg\, t=2,$ $E^{**}_2$ satisfies condition \spd{\Q}{n}. The
statement for higher $r$ follows from Prop \ref{PD-prop}(i) once
we show that $(E^{**}_*,d_*)$ satisfies (ZR). To prove it, choose
a fixed point $x_0\in M$ of the action and consider the diagram
\begin{center}
\begin{tabular}{ccccc}
$\{x_0\}\times_{S^1} ES^1$ & $\stackrel{i}{\longrightarrow}$ &
$M\times_{S^1} ES^1$ &
$\stackrel{\pi}{\longrightarrow}$ & $BS^1$ \\
& $\searrow$ & $\downarrow \pi$ & $\swarrow$ &\\
& & $BS^1$ & &\\
\end{tabular}
\end{center}
where $i$ is the natural embedding, and the
skew arrows represent the identity maps.
Let $(\bar{E}_*^{**}, \bar{d}_*)$ denote the spectral sequence of the
map $id: BS^1\to BS^1,$
$$\bar{E}^{kl}_r=\begin{cases}\Q & \text{for $l=0$ and even $k\geq 0$}\\
0 & \text{otherwise.} \\
\end{cases}$$
The horizontal maps of the diagram above induce morphisms of spectral
sequences $\bar{E}^{**}_* \stackrel{\pi^*}{\longrightarrow} E^{**}_*
\stackrel{i^*}{\longrightarrow} \bar{E}^{**}_*.$
Since $i^*\pi^*$ is the identity on $\bar{E}^{**}_*,$ $E^{k0}_r\ne
0$ for all $r\geq 2$ and any even $k\geq 0.$ Hence
$(E^{**}_*,d_*)$ satisfies condition (ZR).
\end{proof}

Similarly, by using \cite[Thm. 5.2]{McCleary} and adopting the above
argument we prove the following.

\begin{proposition}{\bf PD for Leray-Serre spectral sequence}
If $S^1$ acts on a \pdq{n}-space (with respect to singular
cohomology) $M,$ and $M^{S^1}\ne \emptyset,$ then the Leray-Serre
spectral sequence for singular cohomology of the fibration
$$M\to M\times_{S^1} ES^1\to BS^1$$ satisfies \spd{\Q}{n} for every $r\geq 2.$
\end{proposition}

A convenient way of calculating equivariant cohomology of
a smooth closed manifold $M$ with an $S^1$ action on it is by using the
Cartan construction, \cite[\S4]{AB}:
Let $D^{**}$ be a bigraded $\R$-linear space whose $(2k,l)$th summand
is the space of $S^1$-invariant, differential $l$-forms on
$M,$ $D^{2kl}=\oi{l}\subset \Omega^l(M)$ and $D^{2k+1,l}=0$ for $k\geq 0:$

\begin{center}
\setlength{\unitlength}{3947sp}%
\begingroup\makeatletter\ifx\SetFigFont\undefined%
\gdef\SetFigFont#1#2#3#4#5{%
  \reset@font\fontsize{#1}{#2pt}%
  \fontfamily{#3}\fontseries{#4}\fontshape{#5}%
  \selectfont}%
\fi\endgroup%
\begin{picture}(5075,1602)(537,-3061)
\thinlines
\put(4041,-2578)
{\vector( 3,-1){819.600}}
\put(3870,-2874)
{\vector( 0, 1){315}}
\put(5063,-2874)
{\vector( 0, 1){315}}
\put(5063,-2267)
{\vector( 0, 1){315}}
\put(3870,-2259)
{\vector( 0, 1){315}}
\put(4041,-1966)
{\vector( 3,-1){819.600}}
\multiput(4743,-1723)(-8.11250,2.70417){25}{\makebox(1.6667,11.6667){\SetFigFont{5}{6}{\rmdefault}{\mddefault}{\updefault}.}}\put(4743,-1723)
{\vector( 3,-1){0}}
\put(3870,-1667)
{\vector( 0, 1){196}}
\put(5063,-1667)
{\vector( 0, 1){196}}

\put(3489,-3061){\makebox(0,0)[lb]{\smash{\SetFigFont{11}{13.2}{\rmdefault}{\mddefault}{\updefault}
$\oi{0}$
}}}
\put(4388,-3061){\makebox(0,0)[lb]{\smash{\SetFigFont{11}{13.2}{\rmdefault}{\mddefault}{\updefault}
$0$
}}}
\put(4831,-3061){\makebox(0,0)[lb]{\smash{\SetFigFont{11}{13.2}{\rmdefault}{\mddefault}{\updefault}
$\oi{0}$
}}}
\put(3496,-2454){\makebox(0,0)[lb]{\smash{\SetFigFont{11}{13.2}{\rmdefault}{\mddefault}{\updefault}
$\oi{1}$
}}}
\put(3923,-2776){\makebox(0,0)[lb]{\smash{\SetFigFont{11}{13.2}{\rmdefault}{\mddefault}{\updefault}
$\delta$
}}}
\put(4425,-2640){\makebox(0,0)[lb]{\smash{\SetFigFont{11}{13.2}{\rmdefault}{\mddefault}{\updefault}
$i_X$
}}}
\put(4373,-2468){\makebox(0,0)[lb]{\smash{\SetFigFont{11}{13.2}{\rmdefault}{\mddefault}{\updefault}
$0$
}}}
\put(4816,-2461){\makebox(0,0)[lb]{\smash{\SetFigFont{11}{13.2}{\rmdefault}{\mddefault}{\updefault}
$\oi{1}$
}}}
\put(5123,-2776){\makebox(0,0)[lb]{\smash{\SetFigFont{11}{13.2}{\rmdefault}{\mddefault}{\updefault}
$\delta$
}}}
\put(3481,-1853){\makebox(0,0)[lb]{\smash{\SetFigFont{11}{13.2}{\rmdefault}{\mddefault}{\updefault}
$\oi{2}$
}}}
\put(4823,-1868){\makebox(0,0)[lb]{\smash{\SetFigFont{11}{13.2}{\rmdefault}{\mddefault}{\updefault}
$\oi{2}$
}}}
\put(3938,-2161){\makebox(0,0)[lb]{\smash{\SetFigFont{11}{13.2}{\rmdefault}{\mddefault}{\updefault}
$\delta$
}}}
\put(4388,-1869){\makebox(0,0)[lb]{\smash{\SetFigFont{11}{13.2}{\rmdefault}{\mddefault}{\updefault}
$0$
}}}
\put(4425,-2033){\makebox(0,0)[lb]{\smash{\SetFigFont{11}{13.2}{\rmdefault}{\mddefault}{\updefault}
$i_X$
}}}
\put(5138,-2161){\makebox(0,0)[lb]{\smash{\SetFigFont{11}{13.2}{\rmdefault}{\mddefault}{\updefault}
$\delta$
}}}
\put(1379,-2578)
{\vector( 3,-1){819.600}}
\put(1208,-2874)
{\vector( 0, 1){315}}
\put(2401,-2874)
{\vector( 0, 1){315}}
\put(2401,-2267)
{\vector( 0, 1){315}}
\put(1208,-2259)
{\vector( 0, 1){315}}
\put(1379,-1966)
{\vector( 3,-1){819.600}}
\multiput(2081,-1723)(-8.11250,2.70417){25}{\makebox(1.6667,11.6667){\SetFigFont{5}{6}{\rmdefault}{\mddefault}{\updefault}.}}\put(2081,-1723)
{\vector( 3,-1){0}}
\put(1208,-1667)
{\vector( 0, 1){196}}
\put(2401,-1667)
{\vector( 0, 1){196}}
\put(2721,-1966)
{\vector( 3,-1){819.600}}
\put(2725,-2578)
{\vector( 3,-1){819.600}}
\multiput(3463,-1723)(-8.11250,2.70417){25}{\makebox(1.6667,11.6667){\SetFigFont{5}{6}{\rmdefault}{\mddefault}{\updefault}.}}\put(3463,-1723)
{\vector( 3,-1){0}}
\multiput(5600,-1978)(-8.11250,2.70417){25}{\makebox(1.6667,11.6667){\SetFigFont{5}{6}{\rmdefault}{\mddefault}{\updefault}.}}\put(5600,-1978)
{\vector( 3,-1){0}}
\multiput(5600,-2580)(-8.11250,2.70417){25}{\makebox(1.6667,11.6667){\SetFigFont{5}{6}{\rmdefault}{\mddefault}{\updefault}.}}\put(5600,-2580)
{\vector( 3,-1){0}}
\put(827,-3061){\makebox(0,0)[lb]{\smash{\SetFigFont{11}{13.2}{\rmdefault}{\mddefault}{\updefault}
$\oi{0}$
}}}
\put(1726,-3061){\makebox(0,0)[lb]{\smash{\SetFigFont{11}{13.2}{\rmdefault}{\mddefault}{\updefault}
$0$
}}}
\put(2169,-3061){\makebox(0,0)[lb]{\smash{\SetFigFont{11}{13.2}{\rmdefault}{\mddefault}{\updefault}
$\oi{0}$
}}}
\put(834,-2454){\makebox(0,0)[lb]{\smash{\SetFigFont{11}{13.2}{\rmdefault}{\mddefault}{\updefault}
$\oi{1}$
}}}
\put(1261,-2776){\makebox(0,0)[lb]{\smash{\SetFigFont{11}{13.2}{\rmdefault}{\mddefault}{\updefault}
$\delta$
}}}
\put(1763,-2640){\makebox(0,0)[lb]{\smash{\SetFigFont{11}{13.2}{\rmdefault}{\mddefault}{\updefault}
$i_X$
}}}
\put(1711,-2468){\makebox(0,0)[lb]{\smash{\SetFigFont{11}{13.2}{\rmdefault}{\mddefault}{\updefault}
$0$
}}}
\put(2154,-2461){\makebox(0,0)[lb]{\smash{\SetFigFont{11}{13.2}{\rmdefault}{\mddefault}{\updefault}
$\oi{1}$
}}}
\put(2461,-2776){\makebox(0,0)[lb]{\smash{\SetFigFont{11}{13.2}{\rmdefault}{\mddefault}{\updefault}
$\delta$
}}}
\put(819,-1853){\makebox(0,0)[lb]{\smash{\SetFigFont{11}{13.2}{\rmdefault}{\mddefault}{\updefault}
$\oi{2}$
}}}
\put(2161,-1868){\makebox(0,0)[lb]{\smash{\SetFigFont{11}{13.2}{\rmdefault}{\mddefault}{\updefault}
$\oi{2}$
}}}
\put(1276,-2161){\makebox(0,0)[lb]{\smash{\SetFigFont{11}{13.2}{\rmdefault}{\mddefault}{\updefault}
$\delta$
}}}
\put(1726,-1869){\makebox(0,0)[lb]{\smash{\SetFigFont{11}{13.2}{\rmdefault}{\mddefault}{\updefault}
$0$
}}}
\put(1763,-2033){\makebox(0,0)[lb]{\smash{\SetFigFont{11}{13.2}{\rmdefault}{\mddefault}{\updefault}
$i_X$
}}}
\put(2476,-2161){\makebox(0,0)[lb]{\smash{\SetFigFont{11}{13.2}{\rmdefault}{\mddefault}{\updefault}
$\delta$
}}}
\put(3076,-2640){\makebox(0,0)[lb]{\smash{\SetFigFont{11}{13.2}{\rmdefault}{\mddefault}{\updefault}
$i_X$
}}}
\put(3068,-2033){\makebox(0,0)[lb]{\smash{\SetFigFont{11}{13.2}{\rmdefault}{\mddefault}{\updefault}
$i_X$
}}}
\put(3053,-1869){\makebox(0,0)[lb]{\smash{\SetFigFont{11}{13.2}{\rmdefault}{\mddefault}{\updefault}
$0$
}}}
\put(3054,-3061){\makebox(0,0)[lb]{\smash{\SetFigFont{11}{13.2}{\rmdefault}{\mddefault}{\updefault}
$0$
}}}
\put(3047,-2468){\makebox(0,0)[lb]{\smash{\SetFigFont{11}{13.2}{\rmdefault}{\mddefault}{\updefault}
$0$
}}}
\end{picture}

\end{center}

Let $X\in Vect(M)$ be the vector field on $M$ induced by the $S^1$ action
(the infinitesimal action), let $i_X: \oi{k}\to \oi{k-1}$ be the map
$(i_X\omega)(\cdot,...,\cdot)=\omega(X,\cdot,...,\cdot),$ and let
$\delta$ denote the exterior derivative of differential forms on $M.$
Consider a differential $d:D^{**}\to D^{**},$
$$d(\omega)=\begin{cases} \delta \omega -i_X\omega,& \text{for $\omega\in
D^{kl}=\oi{l}$ for even $k;$}\cr
0 & \text{for odd $k.$}\\
\end{cases}$$
It is not difficult to see that the exterior product of forms
induces a multiplicative structure on $(D^{**},d).$ By Theorem
4.13 and the following paragraphs in \cite{AB}, the cohomology of
the total complex of the above complex is isomorphic to
$H_{S^1}^*(M;R).$ One can show that the vertical filtration of
this double complex yields a spectral sequence which satisfies
\spd{\Q}{n}.


\subsection{Spectral sequences for $\Z/p$-actions}
\label{sSS-Zp}

Let $\Z/p$ act on a paracompact connected space $X.$ Consider the
standard fibration $X\to X_{\Z/p}\stackrel{\pi}{\to} B\Z/p,$ where
$X_{Z/p}=X\times_{\Z/p} B\Z/p.$ There are three spectral sequences
associated with the $\Z/p$-action on $X$
involving cohomology of $X$ with coefficients in a ring $R$:\\
{\bf (Leray)} The Leray spectral sequence of $\pi,$ compare \cite[IV.6]
{Bredon-sheaf}, \cite[5.8.6]{Weibel};\\
{\bf (Serre)} The Leray-Serre spectral sequence of the fibration $\pi$
(defined for singular cohomology theory)
\cite[Ch. 5]{McCleary}, \cite[XIII.7]{Whitehead}.\\
{\bf (Swan)} A spectral sequence defined as follows: Let
$(C^*,\delta)$ be the cochain complex of $X$ for sheaf (or
Alexander-Spanier) cohomology with coefficients in $R.$ There is a
natural $\Z/p$-action on $C^*.$ Let $$\to P^2
\stackrel{\delta'}{\to} P^1 \stackrel{\delta'}{\to} P^0\to R$$ be
a projective resolution of the $R[\Z/p]$-module $R$ with the
trivial $\Z/p$-action. Then $D^{**}=Hom_{R[\Z/p]}(P^*,C^*)$ is a
double complex with the differential $\delta_v=\delta: D^{kl}\to
D^{k,l+1}$ and the differential $\delta_h: D^{kl}\to D^{k+1,l}$
dual to $\delta'.$ We consider the ``first'' spectral sequence
associated with $(D^{**},\delta_h,\delta_v)$ and for the purpose
of this paper we will call it the Swan spectral sequence. (In
\cite{Swan-new} a similar construction based on complete
projective resolutions is considered). A version of Swan spectral
sequence can be constructed for singular cochains of $X$ and for
cellular cochains if $X$ is a CW-complex.

The total complex of $(D^{**},\delta_h,\delta_v)$ is
$$D^s=\bigoplus\limits_{k+l=s} D^{kl},$$
$$d(\alpha)=d_h(\alpha)+(-1)^kd_v(\alpha),$$
where $\alpha\in D^{kl}.$ (Note that $d\circ d=0.$) The ``first''
spectral sequence is the one induced by the vertical filtration of $D^{**}.$

\begin{proposition}\label{spectral_seq}
If $(E^{**}_*,d_*)$ is any of the three spectral sequences
defined above for cohomology with (constant) coefficients in $R$ then
\begin{enumerate}
\item $E_2^{kl}=H^k(\Z/p, H^l(X;R)),$ where $g\in \Z/p$ acts on
$H^l(X;R)$ by the automorphism induced by $g^{-1}:X\to X.$
 \item There is a multiplicative structure $\cdot$ on $(E^{**}_*, d_*)$
such that $\alpha\cdot \beta=(-1)^{k'l}\alpha\cup \beta$ for
$\alpha\in E^{kl}_2,$ $\beta\in E^{k'l'}_2.$ (Note that $H^k(\Z/p,
H^l(X;R)) \times H^{k'}(\Z/p, H^{l'}(X;R))\stackrel{\cup}{\to}
H^{k+k'}(\Z/p, H^{l+l'}(X;R))$ is a well defined cup-product on
cohomology groups with non-constant coefficients.)
 \item If $E^{**}_*$ is the Leray or Leray-Serre spectral sequence then
$E^{**}_*$ converges to $H^*(X_{\Z/p};R).$
\end{enumerate}
\end{proposition}

We do not know if the Swan spectral sequence converges to $H^*(X_{\Z/p};R).$
We do not know either under what conditions the above three
spectral sequences are isomorphic.

{\bf Proof of Proposition \ref{spectral_seq}:}\\
({\bf Leray}) (i) Since all cohomology groups are considered with
supports in closed sets, by \cite[Theorem IV.6.1]{Bredon-sheaf} we
have $E^{kl}_2=H^k(B\Z/p,{\mathcal H}^l(\pi;R)),$ where $B\Z/p$ is
locally contractible and ${\mathcal H}^l(\pi;R)$ is the Leray
sheaf of $\pi.$  By the remark following Theorem IV.8.2 in
\cite{Bredon-sheaf}, ${\mathcal H}^l(\pi;R)$ is locally constant
on $B\Z/p,$ and by careful retracing the relevant definitions, we
see that that Leray sheaf is given by the $\Z/p$-action on
$H^*(X;R)$ described above.
 (ii) follows from
\cite[IV.6.5]{Bredon-sheaf}. (iii) follows from \cite[Theorem
IV.6.1]{Bredon-sheaf}.

\noindent ({\bf Leray-Serre}) The statement for Leray-Serre spectral
sequence follows from
\cite[Thm 5.2]{McCleary}. Compare also \cite[XIII.8.10]{Whitehead}.

\noindent ({\bf Swan}) We have
$E_1^{kl}=Hom_{R[\Z/p]}(P_k,H^l(X;R)),$ where $\Z/p$ acts on
$H^l(X;R)$ as in Proposition \ref{spectral_seq}(i). Therefore,
$E_2^{kl}=H^k(\Z/p,H^l(X;R)).$ If $\Delta: P^*\to P^*\otimes P^*$
is a diagonal approximation of $(P^*,\delta')$ then the cup
product $D^{kl}\otimes D^{k'l'} \stackrel{\cup}{\to}
D^{k+k',l+l'}$ is defined for any $\alpha\in D^{kl}, \beta\in
D^{k'l'}$ by
$$P^{k+k'}\stackrel{\Delta_{kk'}}{\longrightarrow} P^k\otimes
P^{k'} \stackrel{\alpha\otimes \beta}{\longrightarrow} C^l\otimes
C^{l'}\stackrel{\cup}{\to} C^{l+l'}.$$ It has the following
properties:
$$d_h(\alpha\cup \beta)=d_h(\alpha)\cup \beta+ (-1)^k \alpha\cup d_h(\beta),$$
$$d_v(\alpha\cup \beta)=d_v(\alpha)\cup \beta+ (-1)^l \alpha\cup
d_v(\beta),$$ for $\alpha\in D^{kl}, \beta\in D^{k'l'}.$
Let $\alpha \cdot\beta $ be a new product on $D^{**}$ equal to
$(-1)^{k'l}\alpha\cup\beta,$ for $\alpha, \beta$ as above.
The following lemma, whose proof is left to the reader, implies
that $\cdot$ defines a multiplicative structure on $E^{**}_*.$

\begin{proposition}
$$d(\alpha\cdot \beta)=d(\alpha)\cdot \beta+ (-1)^{deg{\alpha}}
\alpha\cdot d(\beta),$$ where $deg(\alpha)=k+l.$
\end{proposition}
This completes the proof of Proposition \ref{spectral_seq}. \Box

Let $(P_*,\delta')$ be the standard resolution of $R$ by free
$R[\Z/p]$-modules: let $P_k=R[\Z/p]$ for all $k\geq 0$ and
let $\delta': P_k\to P_{k-1}$ be
$$\delta'(\alpha)=\begin{cases}(t-1)\cdot \alpha& \text{for $k$ odd,}\\
N\cdot \alpha & \text{for $k$ even}\\
\end{cases}$$
where $t-1, N=1+t+...+t^{p-1}$ are
elements of $R[\Z/p]=R[t]/(t^p-1).$
Let $$\Delta:(P_*,\delta')\to
(P_*,\delta')\otimes (P_*,\delta')$$ be the diagonal approximation
whose
$(k,l)$-component $\Delta_{kl}:P_{k+l}\to P_k\otimes P_l$ is given by
$$\Delta_{kl}(1)=\begin{cases}1\otimes 1 & \text{if $k$ even}\\
1\otimes t & \text{if $k$ odd, $l$ even}\\
\sum_{0\leq i<j\leq p-1} t^i\otimes t^j & \text{if $k,l$ odd,}\\
\end{cases}$$
cf. \cite[Ex. V.1]{Brown}.
Therefore, after identifying
$D^{kl}=Hom_{R[\Z/p]}(P_{k},H^{l}(X;R))$ with $H^l(X;R)$ we have
\begin{equation}\label{product}
\alpha\cup \beta=\begin{cases}\alpha\cup_X \beta & \text{if $k$
is even,} \\
\alpha\cup_X t\beta & \text{if $k$ is odd, $k'$ is even}\\
\sum_{0\leq i<j\leq p-1}t^i\alpha \cup_X t^j\beta& \text{if $k,k'$
are odd,}\\
\end{cases}
\end{equation}
where $\cup_X$ denotes the cup product on $H^*(X;R)$ and $\cup$
denotes the product on $D^{**}$ defined before. (Recall that
$\alpha\cdot \beta=(-1)^{k'l}\alpha\cup \beta,$ for $\alpha\in
D^{kl}, \beta\in D^{k'l'}.$)

\begin{lemma} \label{ZR_for_Zp}
If $R$ is $\Z$ or $\F_p$ and the $\Z/p$-action on $X$ has a fixed
point then all three spectral sequences considered above satisfy
condition (ZR) for $\K=\F_p.$
\end{lemma}

\begin{proof}
Let $E^{**}_*$ be the Leray or Leray-Serre or Swan spectral
sequence associated with the $\Z/p$-action on $X,$ let $x_0\in
X^{\Z/p}$ and let $\bar{E}^{**}_*$ be the corresponding spectral
sequence associated with the trivial $\Z/p$-action on $\{x_0\}.$
The $\Z/p$-equivariant maps: $\{x_0\}\hookrightarrow X$ and $X\to
\{x_0\}$ induce maps
\begin{equation}\label{composition}
(\bar{E}^{**}_r, \bar{d}_r)\to (E^{**}_r, d_r)\to
(\bar{E}^{**}_r, \bar{d}_r),
\end{equation}
whose composition is the identity on $\bar{E}^{**}_r$ for $r\geq 1.$
Since $X$ is assumed connected and $R=\Z$ or $\F_p,$
$E_2^{k0}=H^k(\Z/p,R)$ is either $0$ or $\F_p.$
Hence, if (ZR) is not satisfied then $E^{k0}_2=\F_p$ and
$E^{k0}_\infty=0$ for some
$k.$ This implies that $\bar{E}^{k0}_2=\F_p,$ and since
(\ref{composition}) is the identity map,
$\bar{E}^{k0}_\infty=0.$ This leads to contradiction since
$\bar{E}^{**}_2$ has only one non zero row and $\bar{E}^{**}_\infty=
\bar{E}^{**}_2.$
\end{proof}


\subsection{\P duality for spectral sequences for $\Z/p$-actions}
\label{sPD-Zp}


Let $\Z/p$ act on a \pdfp{n}-space $M$ with a fixed point and let
$(E^{**}_*,d_*)$ be either the Leray or Swan spectral sequence
associated with that action with coefficients in $\F_p.$

\begin{proposition} \label{weakPD}
$(E^{**}_r,d_r)$ satisfies \swpd{\F_p}{n} for $r\geq 2.$
\footnote{V. Puppe pointed to us that a similar result is hidden
in the proof of the main theorem of \cite{Bredon-paper}.}
\end{proposition}

\begin{proof} Let $k,k'$ be of different parity, and let $0\leq l\leq n,$
$l'=n-l.$ Since $\Z/p\subset \Q/\Z,$ and $$H^l(M;\F_p)=
Hom(H^{l'}(M;\F_p),\Z/p)=Hom(H^{l'}(M;\F_p),\Q/\Z)$$ as
$\F_p[\Z/p]$-modules, the duality theorem for Tate cohomology,
\cite[Cor. VI.7.3]{Brown}, implies that
$$H^k(\Z/p,H^l(M;\F_p))\times H^{k'}(\Z/p,H^{l'}(M;\F_p))\to
H^{k+k'}(\Z/p,H^n(M;\F_p))=\F_p$$ is
non-degenerate\footnote{Recall that $k$ and $k'$ have different
parity.}. Therefore, by Proposition \ref{spectral_seq}, $E^{**}_2$
satisfies the weak \P duality. Now the proposition follows from
Lemma \ref{ZR_for_Zp} and Proposition \ref{PD-prop}(ii).
\end{proof}

In order to say more about the multiplicative properties of
$(E^{**}_*, d_*)$ we need to assume that the $\Z/p$-action on $M$
is nice, $H^*(M;\F_p)=T^*\oplus F^*.$ Now
$E^{kl}_2=H^k(\Z/p,H^l(M;\F_p))=T^l$ for $k>0.$ Since the $\Z/p$
action on $T^*$ is trivial, by (\ref{product}) and Proposition
\ref{spectral_seq} for $p\ne 2$ we have
\begin{equation}
\alpha\cdot \beta=\begin{cases}(-1)^{k'l}\alpha\cup\beta &
\text{if $k$ or $k'$ is even,} \\
0& \text{if $k,k'$ are odd}\\
\end{cases}
\label{E_2-product}
\end{equation}
for $\alpha\in E^{kl}_2,\beta\in E^{k'l'}_2.$ (Recall that
$\alpha\cdot \beta=(-1)^{k'l}\alpha\cup\beta.$ Furthermore, for
$k,k'$ odd we have $\alpha\cup\beta=\sum_{0\leq i<j\leq p-1}\alpha
\cup_X \beta=0$ since ${p \choose 2} \equiv 0$ mod $p.$)

\begin{lemma}\label{T-nondeg}
If a $\Z/p$-action on $M$ is nice and $k$ or $k'$ is even then the
product $\cdot$ given by (\ref{E_2-product}) is non-degenerate.
\end{lemma}

\begin{proof}
Since for $k\not\equiv k'$ mod $2$ this follows from Proposition
\ref{weakPD}, we can assume that $k,k'$ are even. Let
$H^l(M;\F_p)=T^l\oplus F^l,$ $H^{l'}(M;\F_p)=T^{l'}\oplus F^{l'},$
and let $\alpha_1,...,\alpha_s$ be generators of the summands of
$F^l=\F_p[\Z/p]\oplus ...\oplus \F_p[\Z/p].$ Let $F^l_i$ be the
$\F_p$-vector subspace of $F^l$ generated by elements
$(t-1)^i\alpha_j$ for $j=1,...,s.$ Note that $F^l=F_0^l\oplus
...\oplus F_{p-1}^l.$ Similarly we decompose $F^{l'}$ into
$F_0^{l'}\oplus ...\oplus F_{p-1}^{l'}.$ Since $M$ is a
\pdfp{n}-space, the matrix representing the product $$(T^l\oplus
F^l_{p-1})\times (T^{l'}\oplus F^{l'}_0\oplus ...\oplus
F^{l'}_{p-1})\stackrel{\cup}{\to} \F_p$$ is of maximal rank,
$dim_{\F_p} T^l+s.$ All columns of this matrix corresponding to
spaces $F^{l'}_i,$ for $i>0,$ are $0.$ Indeed, if $\beta\in
T^l\oplus F^l_{p-1}=(H^l(M;\F_p))^{\Z/p}$ and $\beta'\in F^{l'}_i$
for $i>0$ then there exists $\beta''\in F^{l'}_{i-1}$ such that
$\beta'=(t-1)\beta''.$ Since $$\beta\cup \beta''= t(\beta\cup
\beta'')=\beta\cup t\beta'',$$ we have $\beta\cup
(t-1)\beta''=\beta\cup \beta'=0.$ Therefore, the matrix of the cup
product on
$$T^l\oplus F^l_{p-1}\times T^{l'}\oplus F^{l'}_0$$
is non-degenerate. By an argument similar to the above,
$\beta\cup \beta'=0$ for any $\beta\in F^l_{p-1}$ and $\beta'\in T^{l'}.$
Hence this matrix has a form
$$\begin{array}{cc} & T^{l'}\ F^{l'}_0\\
\begin{array}{c} T^l\\ F^l_{p-1}\end{array} &
\left(\begin{array}{cc} A & B\\ 0 & C\end{array}\right).
\end{array}$$
Therefore, the matrix $A$ associated with $T^l\times T^{l'}
\stackrel{\cup}{\to} \F_p$ is non-degenerate.
\end{proof}

The above lemma shows that $E^{**}_2$ satisfies \sspd{\F_p}{n} for
nice $\Z/p$-actions. If $p=2$ then it can be shown by induction on
$r$ and by using Proposition \ref{PD-main-prop} that $E^{**}_r$
satisfies \sspd{\F_p}{n} for all $r\geq 2.$ However, we do not
know if $E^{**}_r$ satisfies \sspd{\F_p}{n} for $p\ne 2$ in
general. This problem stems from the fact that the implication of
Proposition \ref{PD-prop}(iii) does not hold for odd $r.$
Therefore, for certain applications it is necessary to assume
condition (\ref{cond}).

\begin{lemma}\label{SSPD}
If condition (\ref{cond}) holds for a given $\Z/p$-action on $M$
then $(E^{**}_r,d_r)$ satisfies \sspd{\F_p}{n} for each $r\geq 2.$
\end{lemma}

\begin{proof} The statement follows from  Proposition \ref{PD-prop}(iii)
and Lemmas \ref{ZR_for_Zp} and \ref{T-nondeg}.
\end{proof}

\begin{proposition} \label{SPD_3mflds}
Condition (\ref{cond}) holds for $n\leq 3.$ Consequently, for any
nice  $\Z/p$-action on a \pdfp{n}-space for $n\leq 3$ each term of
the induced Leray spectral sequence satisfies \sspd{\F_p}{n} for
all $r.$
\end{proposition}

\begin{proof} For $n=1$ the statement is obvious. For $n=2$
the statement is a consequence of condition (ZR), cf. Lemma
\ref{ZR_for_Zp}. Therefore, assume that $n=3.$ Since $d_r=0$ for
$r\geq 5,$ it suffices to show that $d^{kl}_3=0$ for $k\geq n.$
For $l=2$ it follows from Lemma \ref{ZR_for_Zp}. Hence assume that
$l=3$ and that $d^{kl}_3(\omega)=\alpha \ne 0$ for some $\omega\in
E^{k3}_3=\F_p.$ Since $\alpha\in E^{k+3,1}_3,$ by the weak \P
duality there exists $\beta\in E^{k+2,2}_3$ such that $\alpha\cdot
\beta\ne 0.$ By Lemma \ref{ZR_for_Zp}, $d_3(\beta)=0$ and hence we
get a contradiction:
$$0=d_3(\omega\beta)=d_3(\omega)\cdot \beta+\omega\cdot d_3(\beta)=
\alpha\cdot\beta\ne 0.$$
\end{proof}

The next result concerns \P duality for spectral sequences
with integral coefficients.

\begin{proposition}\label{PD-Z}
If $\Z/p$ acts nicely on a \pdfp{n}-space $M,$ with no $p$-torsion
in $H^*(M;\Z)$ and if $M^{\Z/p}\ne \emptyset$ then the Leray and
the Swan spectral sequences for that action and for $R=\Z$ satisfy
\P duality, \spd{\F_p}{n}, for all $r\geq 2.$
\end{proposition}

\begin{proof} The statement for $r=2$ follows from Lemmas \ref{nice}
and \ref{T-nondeg}. For $r> 2$ the statement follows from Lemma
\ref{ZR_for_Zp} and Proposition \ref{PD-prop}(i).
\end{proof}


\section{Proofs of the main results}
\label{sProofs}


\begin{lemma}\label{mod_4}
Let $(E^{**}_*,d_*)$ be a spectral sequence whose terms
$(E^{**}_r,d_r)$ for $r\geq 2$ are vector spaces over a field
$\K,$ $char\, \K\ne 2,$ and satisfy either \spd{\K}{n} or
\sspd{\K}{n}. In the latter case we assume that $d_r=0$ for odd
$r.$ If
 \begin{enumerate}
 \item $n$ is even, or
 \item $E^{kl}_2=0$ for all even $l,$
$0<l\leq \frac{1}{2}(n-1),$ and for all sufficiently large $k$
\end{enumerate}
then
$$\sum_l dim_{\K}\, E^{kl}_\infty\equiv
\sum_l dim_{\K}\, E^{kl}_2\quad \mod 4,$$
for all sufficiently large $k.$
\end{lemma}

\begin{proof}
It is enough to prove that
$$\sum_l dim\, E^{kl}_{r+1}\equiv \sum_l dim\, E^{kl}_r\quad mod\, 4,$$
for all $r\geq 2$ and sufficiently large $k.$ Since
$E^{**}_r=E^{**}_{r+1}$ for $r$ odd, assume that $r$ is even. By
Lemma \ref{prop_of_PD}, $rank\, d^{kl}_r=rank\, d^{k-r,l}_r,$ for
all sufficiently large $k$. Therefore,
$$\sum_{l} dim\, E^{kl}_{r+1}=
\sum_{l} dim\, Ker\, d^{kl}_r-
\sum_{l} dim\, Im\, d^{k-r,l+r-1}_r=$$
$$\sum_{l} dim\, E^{kl}_r -
\sum_{l} rank\, d^{kl}_r-\sum_{l} rank\, d^{k-r,l+r-1}_r=$$
$$\sum_{l} dim\, E^{kl}_r -
2\sum_{l} rank\, d^{kl}_r.$$
Therefore, we need to prove that
$$\sum_{l} rank\, d^{kl}_r\equiv 0\quad \mod 2.$$
By Lemma \ref{prop_of_PD},
$$\sum_{l} rank\, d^{kl}_r=2\cdot \hspace*{-.1in}
\sum_{l<n-l+r-1}\hspace*{-.1in}  rank\, d^{kl}_r +
\begin{cases} rank\, d^{kl_0}_r & \text{if $l_0=n-l_0+r-1$}\\
0 & \text{if there is no such $l_0.$}\\
\end{cases}$$
For $n$ even, there is no such $l_0$ and the proof is completed.
Hence, assume that $n$ is odd. If $l_0$ is odd then
$l_0-r+1$ is even and $l_0-r+1\leq \frac{1}{2}(n-1).$ Hence
$E^{k+r,l_0-r+1}_r=0$ by the assumption of the lemma, and therefore
$d^{kl_0}_r: E^{kl_0}_r\to E^{k+r,l_0-r+1}_r$ is $0.$
Therefore assume that $l_0$ is even. Since $d^{kl}_r=0$ for odd $k,$ assume
also that $k$ is even.
Consider the bilinear form
$$\Psi: E_r^{kl_0}\times E_r^{kl_0}\to
E_r^{2k+r,n}=\K,$$
$\Psi(\alpha,\beta)=d_r(\alpha)\cdot \beta.$
We have
$$d_r(\alpha)\cdot \beta + (-1)^{k+l_0} \alpha\cdot d_r(\beta)=
d_r(\alpha\cdot \beta)=0.$$ Since $deg(\alpha)=k+l_0$ is even,
$\alpha\cdot d_r(\beta)= d_r(\beta)\cdot\alpha$ and
$$d_r(\alpha)\cdot \beta + d_r(\beta)\cdot \alpha=0.$$ Therefore
$\Psi$ is skew-symmetric, and it has an even rank. But $rank\,
\Psi=rank\, d^{kl_0}_r,$ since $\alpha\in E^{kl_0}_r,$
$d_r(\beta)\in E^{k+r,l_0-r+1}_r$ and the product\\
$E^{kl_0}_r\times E^{k+r,l_0-r+1}_r \to E^{2k+r,n}$ is
non-degenerate.
\end{proof}


\subsection{Proof of Theorem \ref{torus}}
\label{sTorus-proof}


The following lemma shows that it is sufficient to prove Theorem \ref{torus}
for circle actions.

\begin{lemma} If an action of a torus $T$ on $X$ has FMCOT then
there exists $S^1\subset T$ such that $X^{S^1}=X^T.$
\end{lemma}

\begin{proof}
The condition FMCOT implies that the set $\{(T_x)^0: x\in X\}$ is finite.
Denote its elements different than $T$ by $T_1,...,T_n.$ Consider
$S^1\subset T$ which does not lie inside $T_i$ for any $i.$
Then $S^1\cap T_i$
is finite for $1\leq i\leq n.$ Since each $T_i$ has only countably many
finite extensions in $T,$ the set $S^1\cap \left(\bigcup_{x\in X\setminus X^T}
T_x\right)$ is at most countable. Therefore there exists $t\in S^1$
such that the only points of $X$ fixed by $t$ are the elements of $X^T.$
Hence $X^{S^1}=X^T.$
\end{proof}

Assume now that $T=S^1.$ The proof of Theorem \ref{torus} for $n$
even follows immediately from (\ref{euler-torus}), the lemma
below, and the fact that $M^T$ has finitely many components, each
of which is a \pdq{m}-space, for $m$ even (see Theorem 5.2.1 and
Remark 5.2.4 in \cite{AP}; cf. \cite{CS}).

\begin{lemma}\label{2n} If $n$ is even and $M$ is a \pdq{n}-space,
then\\
$\sum_i b^i(M)\equiv \chi(M)$ mod $4.$
\end{lemma}

\begin{proof}
Since $b^i(M)=b^{n-i}(M),$ the difference between the left and the right
side of the above identity is
$$2\sum_{odd\ i} b^i(M)=4\sum_{odd\ i<n/2} b^i(M)+
2\begin{cases}b^{n/2}(M) & \text{if $n/2$ is odd} \\ 0 & \text{otherwise.}\\
\end{cases}$$
This completes the proof for $n/2$ even. If $n/2$ is odd then
the pairing $$H^{n/2}(M;\Q)\times H^{n/2}(M;\Q) \stackrel{\cup}{\to}
H^n(M;\Q)=\Q$$ is non-degenerate and skew-symmetric. Hence
$b^{n/2}(M)$ is even.
\end{proof}

Assume now that $n$ is odd and that an action of $S^1$ on a \pdq{n}-space
$M$ satisfies all assumptions of Theorem \ref{torus}.

By Proposition 3.10.9 and Corollary 3.10.12 in \cite{AP},
$H^i_{S^1}(M,M^{S^1};\Q)=0$ for $i>cd\, M.$
Therefore, the long exact sequence of the equivariant cohomology groups
for the pair $(M,M^{S^1})$ gives an isomorphism

\begin{equation}
H^s_{S^1}(M;\Q)=H^s_{S^1}(M^{S^1};\Q),
\label{localization}
\end{equation}
for $s>cd\, M.$

$H^*_{S^1}(M^{S^1};\Q)=H^*(M^{S^1};\Q)\otimes H^*(BS^1)$ and
$H^*(BS^1)=\Q[t],$ where $deg\, t=2.$ Therefore,
\begin{equation}
dim_\Q\, H^s_{S^1}(M^{S^1};\Q)+dim_\Q\, H^{s+1}_{S^1}(M^{S^1};\Q)=
\sum_i dim_\Q\, H^i(M^{S^1};\Q).
\label{right}
\end{equation}

The Leray spectral sequence $(E^{**}_*,d_*)$ of the map
$M\times_{S^1} ES^1\to BS^1$ with coefficients in $\Q$ converges
to $H^*_{S^1}(M;\Q).$ By Proposition \ref{PD-S^1},
$(E^{**}_r,d_r)$ satisfies condition \spd{\Q}{n} for all $r\geq
2.$ Since by Lemma \ref{prop_of_PD}(i) the ranks of entries in
$E^{**}_\infty$ are $2$-periodic, we have
$$dim\, H^s_{S^1}(M;\Q)+dim\, H^{s+1}_{S^1}(M;\Q)=$$
$$\hspace*{-.1in}
\sum_{k+l=s\atop k\, even} \hspace*{-.1in} dim\, E^{kl}_{\infty}+
\hspace*{-.1in} \sum_{k+l=s+1\atop k\, even} \hspace*{-.1in} dim\,
E^{kl}_{\infty}= \sum_l dim\, E^{k_0l}_{\infty}$$ for sufficiently
large $s$ and sufficiently large even $k_0.$ By Lemma \ref{mod_4}
the above expression is equal mod $4$ to $\sum_l dim\, E^{k_0l}_2=
\sum_l dim\, H^l(M;\Q).$ Hence $$dim\, H^s_{S^1}(M;\Q)+dim\,
H^{s+1}_{S^1}(M;\Q)\equiv \sum_l dim\, H^l(M;\Q)\quad \mod 4.$$
This equality together with (\ref{localization}) and (\ref{right})
implies Theorem \ref{torus}.\Box

\subsection{Proof of Proposition \ref{p-euler-zp}}
\label{speuler-zp}


The following three lemmas will be needed in the proof of
Proposition \ref{p-euler-zp} and of Theorem \ref{Zp}:

\begin{lemma} \label{quillen} If $\Z/p$ acts on a paracompact space
$X$ of finite cohomological dimension then the embedding
$i:X^{\Z/p}\to X$ induces an isomorphism $i^*:H_{\Z/p}^s(X;A) \to
H_{\Z/p}^s(X^{\Z/p};A)$ for $s>cd\, X,$ where $A$ is an arbitrary
group of (constant) coefficients.
\end{lemma}

\begin{proof} By Proposition 3.10.9 in \cite{AP},
$H^*_{\Z/p}(X,X^{\Z/p};A)\simeq$\\ $H^*(X/(\Z/p),
X^{\Z/p}/(\Z/p);A).$ Since by \cite[Prop. A.11]{Qu},
 $cd\, (X/(\Z/p))\leq cd\, X,$ we have $H^s(X/(\Z/p), X^{\Z/p}/(\Z/p);A)=0$
for $s>cd\, X.$ Now the proposition follows from the long exact
sequence for the equivariant cohomology of the pair $(X,
X^{\Z/p}).$
\end{proof}

\begin{lemma}\label{k-formula}
{\bf K\"unneth formula for sheaf cohomology} If $R$ is a principal
ideal domain, $Y$ is a CW-complex and $X$ is a paracompact space
such that $H^l(X;R)$ is finitely generated $R$-module for each
$l,$ then there exists a split exact sequence\vspace*{.1in}\\
$0\to \bigoplus_{k+l=s} H^k(X;R)\otimes H^l(Y;R) \to H^s(X\times
Y;R)\to$
\begin{flushright}
$\bigoplus_{k+l=s+1} Tor_R(H^k(X;R), H^l(Y;R))\to 0.$
\end{flushright}
\end{lemma}

\begin{proof}
By \cite[Prop A.4]{Hatcher} (cf. \cite[Ex. Ch7 E5]{Spanier}) $Y$
is a locally contractible space. By the remark following Theorem
IV.8.2 in \cite{Bredon-sheaf}, the Leray sheaf of the projection
$\pi: X\times Y\to Y$ is the constant sheaf with the stalk
$H^*(X;R).$ Therefore, the statement of proposition follows from
\cite[Ex. IV.18]{Bredon-sheaf}.
\end{proof}

\begin{lemma}\label{kunneth} If $\Z/p$ acts trivially on $X$ then
$$H^s_{\Z/p}(X;\Z)\cong \bigoplus_{l\equiv s \smod 2}
H^l(X;\F_p)$$ for $s>cd\, X.$
\end{lemma}

\begin{proof}
By Lemma \ref{k-formula}, $H^s_{\Z/p}(X;\Z)$ is isomorphic to
$$\bigoplus_{k+l=s}
H^k(B\Z/p,\Z)\otimes H^l(X;\Z)\oplus \bigoplus_{k+l=s+1}
Tor_{\Z}(H^k(B\Z/p,\Z), H^l(X;\Z)).$$ Since for $k>0$
$H^k(B\Z/p,\Z)$ is either $\F_p$ or $0$ depending if $k$ even or
odd, $$H^s_{\Z/p}(X;\Z)=\bigoplus_{l\equiv s \smod 2}
H^l(X;\Z)\otimes \F_p\oplus \bigoplus_{l\equiv s-1 \smod 2}
Tor_{\Z}(H^l(X;\Z),\F_p),$$ for $s>cd\, X.$ But by the universal
coefficient theorem for cohomology, the right side is isomorphic
to $\bigoplus_{l\equiv s \smod 2} H^l(X;\F_p).$
\end{proof}

For the proof of Proposition \ref{p-euler-zp} we will need the
following version of the notion of Euler characteristic for double
complexes: if $D^{**}$ is a double complex of vector spaces over a
field $\K$ then let $$\chi(D^{**})=\lim_{N\to \infty} \frac{1}{N}
\sum_{0\leq k\leq N \atop l\in \Z} (-1)^{k+l} dim_\K\, D^{kl}$$ if
this limit exists.

\begin{proposition}\label{euler} If $(E^{**}_*,d_*)$ is a spectral sequence
such that (a) $E^{kl}_r$ are vector spaces over a field $\K$ and
$dim_\K\, E^{kl}_r\leq c$ for all $k,l,$ for a certain $c,$ (b)
$E^{*l}_r=0$ for all $l<0$ and $l>n$ for some $n,$ (c)
$\chi(E^{**}_r)$ exists,
then\\
(i) $E^{**}_{r+1}$ satisfies conditions (a),(b),(c) as well, and\\
(ii) $\chi(E^{**}_{r+1})=\chi(E^{**}_r).$
\end{proposition}

\begin{proof}
The only non trivial statement of the proposition is that
$\chi(E^{**}_{r+1})$ exists and it is equal to
$\chi(E^{**}_r).$ Consider the cochain complex $(C^*_{kl,r}, d_r)$, where
$C^i_{kl,r}=E^{k+ir,l-i(r-1)}_r$ for $i\in \Z.$
Note that under the above assumptions the sums
$$\sum_{0\leq k\leq N \atop l\in \Z} (-1)^{k+l} dim_\K\, E^{kl}_r$$
and $$\sum_{0\leq k\leq N \atop 0\leq l\leq r-1} (-1)^{k+l} dim_\K\,
\chi(C^*_{kl,r})$$ differ by a finite number of terms of the form
$(-1)^{k+l}dim_\K E^{kl}_r$, and that the number of such terms
does not depend on $N.$
Since $dim_\K\, E^{kl}_r\leq c,$ the difference between the above two sums
is bounded uniformly in $N$ and hence
\begin{equation}\label{euler1}
\chi(E^{**}_r)=\lim_{N\to \infty}\, \frac{1}{N}
\sum_{0\leq k\leq N \atop 0\leq l\leq r-1} (-1)^{k+l} dim_\K\,
\chi(C^*_{kl,r}).
\end{equation}
Since $\chi(C^*_{kl,r})=\chi(C^*_{kl,r+1}),$
the proof follows from (\ref{euler1}) and the analogous equation for $r+1.$
\end{proof}

Let $p\ne 2$ and let $\Z/p$ act nicely on a space $X$ with no $p$-torsion in
$H^*(X;\Z),$ and let $(E^{**}_*,d_*)$ be the associated Leray spectral
sequence with coefficients in $\Z.$
Since $E^{kl}_2=H^k(\Z/p,H^l(X;\Z)),$ by Lemma \ref{nice},
$\chi(E^{**}_2)=\frac{1}{2}\chi_t(X).$ Therefore, by Proposition \ref{euler},
$\chi(E^{**}_\infty)$ exists and
\begin{equation}\label{chi-E}
\chi(E^{**}_\infty)=\frac{1}{2}\chi_t(X)
\end{equation}
By an argument similar to that used in the proof above,
$$\chi(E^{**}_\infty)=
\lim_{N\to \infty}\, \frac{1}{N} \sum_{0\leq s\leq N} (-1)^s
\sum_k dim_{\F_p} E^{k,s-k}_\infty.$$ By Lemmas \ref{quillen} and
\ref{kunneth},
$$\sum_k dim_{\F_p} E^{k,s-k}_\infty \cong
Gr\, H^s_{\Z/p}(X;\Z)\cong Gr\, H^s_{\Z/p}(X^{\Z/p};\Z) \cong
\hspace{-.1in} \bigoplus_{l\equiv s \smod 2}\hspace*{-0.1in}
H^l(X^{\Z/p};\F_p)$$ for $s>cd\, X.$ Hence
$$\chi(E^{**}_\infty)=\frac{1}{2}\chi(H^*(X^{\Z/p};\F_p)).$$ Since
$H^2(\Z/p, H^l(X^{\Z/p};\F_p))=H^l(X^{\Z/p};\F_p),$
$$\chi(E^{**}_\infty)=\frac{1}{2}\chi_t(X^{\Z/p}).$$
Now, by (\ref{chi-E}), the proof is completed.

\subsection{Proof of Theorem \ref{Zp}}
\label{sZp-Z-proof}


Let $\Z/p$ act on a \pdfp{n} space $M$ in such a way that the
assumptions of Theorem \ref{Zp} are satisfied. If $n$ is even,
than Theorem \ref{Zp} can be given a proof analogous to that for
$S^1$-actions. Indeed, by Theorem 5.2.1 and Remark 5.2.4 in
\cite{AP} (cf. \cite{CS}), $M^{\Z/p}$ has finitely many
components, each of which is a \pdfp{m}-space, for $m$ even.
Therefore, by Proposition \ref{p-euler-zp} it is enough to prove
that $\sum_i t^i(M)\equiv \chi_t(M)$ mod $4.$ As in the proof of
Lemma \ref{2n}, it is sufficient to show that $t^{n/2}(M)$ is even
if $n/2$ is odd. By Lemma \ref{T-nondeg}, the cup product on
$H^2(\Z/p,H^{n/2}(M;\F_p))$ is non-degenerate and, since it is
skew-symmetric, $t^{n/2}(M)=$\\ $dim_{\F_p}\,
H^2(\Z/p,H^{n/2}(M;\F_p))$ is even.

Assume now that $n$ is odd and $M^{\Z/p}\ne \emptyset.$ Consider
the Leray spectral sequence, $(E^{**}_*,d_*)$, with coefficients
in $\Z$ associated with the $\Z/p$-action on $M.$ Since
$(E^{**}_*,d_*)$ converges to $H^*_{\Z/p}(M;\Z),$ there is a
filtration of $H^*_{\Z/p}(M;\Z)$ such that
$$Gr\, H^s_{\Z/p}(M;\Z)=\bigoplus_i F^i H^s_{\Z/p}(M;\Z)/
F^{i+1} H^s_{\Z/p}(M;\Z)=\bigoplus_{k+l=s} E^{kl}_\infty.$$ Since
$H^l(M;\Z)$ is finitely generated, $E_2^{kl}=H^k(\Z/p,H^l(M;\Z))$
is a finite dimensional vector space over $\F_p$ for $k>0.$

\begin{corollary}\label{p-group} If $s>cd\, X$ then
$H^s_{\Z/p}(M;\Z)$ is a finite $p$-group and\\
$Gr\, H^s_{\Z/p}(M;\Z)$ is a finite dimensional vector space over
$\F_p.$
\end{corollary}

By Proposition \ref{PD-Z}, $(E^{**}_*,d_*)$ satisfies \P
duality $P_{\F_p}(n).$ Therefore, by Lemma \ref{prop_of_PD}(i), we
have
$$dim_{\F_p}\, Gr\, H^s_{\Z/p}(M;\Z)+dim_{\F_p}\, Gr\, H^{s+1}_{\Z/p}(M;\Z)=
\sum_l dim_{\F_p} E^{kl}_\infty,$$ for sufficiently large $s$ and
$k,$ $k$ even. By the assumptions of Theorem \ref{Zp} and by Lemma
\ref{nice}, $E^{kl}_2=0$ for all even $l,$ $0< l\leq
\frac{1}{2}(n-1)$ and for $k>0.$ Therefore, by Proposition
\ref{PD-Z}, the assumptions of Lemma \ref{mod_4} are satisfied and
the sum above equals to $\sum dim_{\F_p}\, E^{kl}_2$ mod $4.$
Therefore, by Lemma \ref{nice},
\begin{equation}
dim_{\F_p}\, Gr\, H^s_{\Z/p}(M;\Z)+dim_{\F_p}\, Gr\,
H^{s+1}_{\Z/p}(M;\Z) \equiv \sum_l t^l(M) \mod 4, \label{one-side}
\end{equation}
for sufficiently large $s.$

By Lemmas \ref{quillen} and \ref{kunneth}, the left hand side of
(\ref{one-side}) is equal to\\ $\sum_l dim_{\F_p}\,
H^l(M^{\Z/p};\F_p).$ Hence by (\ref{t_for_trivial}),
$$\sum_l t^l(M^{\Z/p})=\sum_l dim_{\F_p}\, H^l(M^{\Z/p};\F_p)\equiv
\sum_l t^l(M) \mod 4.$$ \Box

\subsection{Proof of Theorem \ref{Zp-Fp}}
\label{sZp-Fp-proof}


Let $\Z/p$ act on $M$ such that all assumptions of Theorem
\ref{Zp-Fp} are satisfied. The proof of Theorem \ref{Zp-Fp} is
analogous to the proof of Theorem \ref{Zp}, except that $\F_p$ is
the ring of coefficients this time. Let $(E^{**}_*,d_*)$ be the
Leray spectral sequence with coefficients in $\F_p$ associated
with the $\Z/p$-action on $M.$ It converges to
$H^s_{\Z/p}(M;\F_p).$ Since $H^s_{\Z/p}(M;\F_p)$ is a vector space
over $\F_p,$ $Gr\, H^s_{\Z/p}(M;\F_p)\cong H^s_{\Z/p}(M;\F_p)$ for
any filtration of $H^s_{\Z/p}(M;\F_p),$ and hence
$$H^s_{\Z/p}(M;\F_p)\cong \sum_{k+l=s} E^{kl}_\infty.$$
By Lemmas \ref{SSPD} and \ref{prop_of_PD}(ii) the dimensions of
$E^{kl}_r$ over $\F_p$ do not depend on $k$ for large $k$.
Therefore, by Lemma \ref{mod_4}
\begin{equation}\label{zp-fp-1}
\begin{array}{l}
dim_{\F_p}\, H^s_{\Z/p}(M;\F_p)=\sum_l dim_{\F_p}\, E^{k_0l}_\infty\equiv \\
\hspace*{1in} \sum_l dim_{\F_p}\, E^{k_0l}_2=\sum_l t^l(M)\quad \mod 4,
\end{array}
\end{equation}
for $s>cd\, M$ and large $k_0.$

On the other hand, by Lemma \ref{k-formula}
$$H^s_{\Z/p}(M^{\Z/p};\F_p)\cong \bigoplus_{k+l=s}
H^k(B\Z/p,\F_p)\otimes H^l(M^{\Z/p};\F_p)=\bigoplus_{l} H^l(M^{\Z/p};\F_p).$$
Therefore,
\begin{equation}\label{zp-fp-2}
dim_{\F_p}\, H^s_{\Z/p}(M^{\Z/p};\F_p)=\sum_l dim_{\F_p}\, H^l(M^{\Z/p};\F_p)
=\sum_l dim_{\F_p}\, t^l(M^{\Z/p}).
\end{equation}
By Lemma \ref{quillen}, the left sides of (\ref{zp-fp-1}) and
(\ref{zp-fp-2}) are equal, and hence the proof is completed. \Box

\centerline{Department of Mathematics, 244 Mathematics Bldg.}
\centerline{State University of New York at Buffalo, Buffalo NY
14260} \centerline{asikora@buffalo.edu}\vspace*{0.1in}

\centerline{and}\vspace*{0.1in}

\centerline{Institute for Advanced Study, School of Mathematics}
\centerline{1 Einstein Dr., Princeton NJ 08540}

\begin{thebibliography}{99}
\bibitem[AD] {Adem} A. Adem, J. Davis, Topological Transformation
Groups, preprint, http://www.arxiv.org/abs/math.AT/9706228
\bibitem [AHP] {AHP} C. Allday, B. Hanke, V. Puppe, \P duality in
P. A. Smith theory, {\em Proc. of AMS,} to appear.
\bibitem [AB] {AB} M. F. Atiyah, R. Bott, The Moment Map And Equivariant
Cohomology, Topology, {\bf 23} (1984) no. 1, 1--28.
\bibitem [AP] {AP} C. Allday, V. Puppe, Cohomological Methods in
Transformation Groups, Cambridge studies in advanced mathematics 32,
Cambridge Univ. Press, 1993.
\bibitem[BT]{Bott-Tu} R. Bott, L. W. Tu, Differential forms in algebraic
topology, Graduate texts in mathematics 82, Springer-Verlag, 1982.
\bibitem[Br1]{Bredon-groups} G. E. Bredon, Introduction to Compact
Transformation Groups, Academic press 1972.
\bibitem[Br2]{Bredon-sheaf} G. E. Bredon, Sheaf Theory, McGraw-Hill
series in higher mathematics, 1967.
\bibitem[Br3]{Bredon-paper} G. E. Bredon, Fixed point sets of actions on
\P duality spaces, {\em Topology} {\bf 12} (1973), 159-175.
\bibitem[Bro]{Brown} K. S. Brown, Cohomology of groups,
Graduate texts in mathematics 87, Springer Verlag, 1982.
\bibitem[Bry]{Bryan} J. Bryan, Undergraduate Thesis, Stanford Univ., 1989.
\bibitem[Bu]{Bourbaki} N. Bourbaki, General Topology, Part 2,
Addison-Wesley 1966.
\bibitem[CS]{CS} T. Chang, T. Skjelbred, Group actions on \P
duality spaces, {\em Bull. AMS} {\bf 78} (1972), 1024--1026.
\bibitem[CR]{Curtis} C. Curtis, I. Reiner, Representation Theory of
Finite Groups, Interscience Publishers, 1962.
\bibitem[Hat]{Hatcher} A. Hatcher, Algebraic Topology,
Cambridge Univ. Press, 2002. See also\\
http://www.math.cornell.edu/$\sim$hatcher
\bibitem[Han]{Hanke} B. Hanke, Realization of Equivariant Chain Complexes,
{\em Topology,} to appear.
\bibitem [McC]{McCleary} J. McCleary, User's Guide to Spectral Sequences,
Publish or Perish Inc. 1985.
\bibitem [OR]{OR} P. Orlik and F. Raymond, Actions of $SO(2)$ on
$3$-manifolds, in Proceedings of the Conference on Transformation
Groups, Springer Verlag, 1968, 297--318.
\bibitem [PS] {PS} J. H. Przytycki and M. Sokolov, Surgeries on periodic
links and homology of periodic $3$-manifolds, {\em Math. Proc.
Camb. Phil. Soc.} {\bf 131}(2) (2001) 295--307. See also
``Some corrections written after the paper was published'' at\\
http://gwis2.circ.gwu.edu/$\sim$przytyck/publications.
\bibitem[Qu]{Qu} D. Quillen, The spectrum of an equivariant
cohomology ring I, {\em Ann. of Math.} {\bf 94} (3) 1971, 549--572.
\bibitem [Ray]{Ra} F. Raymond, Classification of the actions of the
circle on $3$-manifolds, {\em Trans. Amer. Math. Soc.} {\bf 131} (1968),
51--78.
\bibitem[Sp]{Spanier} E. H. Spanier, Algebraic Topology, Springer 1966.
\bibitem[Sw1]{Swan-new} R. G. Swan, A new method in the fixed point theory,
Comment. Math. Helv. {\bf 34} (1960), 1--16.
\bibitem[Sw2]{Swan} R. G. Swan, The Theory of Sheaves, Chicago Lectures
in Mathematics, Univ. of Chicago Press, 1964.
\bibitem[We]{Weibel} C. A. Weibel, An introduction to homological
algebra, Cambridge studies in advanced mathematics 36, Cambridge
Univ. Press, 1994.
\bibitem[Wh]{Whitehead} G. W. Whitehead, Elements of Homotopy Theory,
Graduate Texts in Mathematics, Springer-Verlag 1978.\vspace*{1in}
\end{thebibliography}
\end{document}